\documentclass{emsprocart}

\contact[bartelsa@math.uni-muenster.de]{Arthur Bartels, Mathematisches
  Institut, Westf\"alische Wilhelms-Universit\"at M\"unster,
  Einsteinstr.~62, D-48149 M\"unster, Germany}
\contact[lueck@math.uni-muenster.de]{Wolfgang L\"uck, Mathematisches
  Institut, Westf\"alische Wilhelms-Universit\"at M\"unster,
  Einsteinstr.~62, D-48149 M\"unster, Germany}
\contact[holger.reich@googlemail.com]{Holger Reich, Mathematisches
  Institut, Heinrich-Heine-Universit\"at D\"usseldorf,
  Universit\"atsstr.~1, D-40225 D\"usseldorf, Germany}

\title[On the Farrell-Jones Conjecture and its applications]{On the
  Farrell-Jones Conjecture and its applications}

\author[Arthur Bartels, Wolfgang L\"uck and Holger Reich]{Arthur
  Bartels, Wolfgang L\"uck and Holger Reich}

\newtheorem{theorem}{Theorem}[section]
\newtheorem{corollary}[theorem]{Corollary}
\newtheorem{lemma}[theorem]{Lemma}

\newtheorem{conjecture}[theorem]{Conjecture}

\theoremstyle{definition} \newtheorem{definition}[theorem]{Definition}
\newtheorem{remark}[theorem]{Remark}

\numberwithin{equation}{section}

\usepackage{enumerate,amssymb,amscd}

\usepackage[arrow,curve,matrix,tips,2cell]{xy} \SelectTips{eu}{10}
\UseTips \UseAllTwocells

\newcounter{commentcounter}

 \DeclareMathOperator{\aut}{aut}

 \DeclareMathOperator{\con}{con}

\DeclareMathOperator{\End}{End} 
 \DeclareMathOperator{\id}{id}
\DeclareMathOperator{\ind}{ind} 
 
 \DeclareMathOperator{\pr}{pr}
 
\DeclareMathOperator{\tr}{tr}

\DeclareMathOperator{\im}{im}

\DeclareMathOperator*{\colim}{colim}

\newcommand{\Fin}{{\mathcal{F}\text{in}}}
\newcommand{\FCyc}{{\mathcal{FC}\text{yc}}}
\newcommand{\VCyc}{{\mathcal{VC}\text{yc}}}

\newcommand{\IC}{\mathbb{C}} 
 \newcommand{\IF}{\mathbb{F}}

 \newcommand{\IQ}{\mathbb{Q}}
\newcommand{\IR}{\mathbb{R}}

\newcommand{\IZ}{\mathbb{Z}}

\newcommand{\calc}{\mathcal{C}} 
 \newcommand{\calf}{\mathcal{F}}
\newcommand{\calfj}{\mathcal{F}\!\mathcal{J}} \newcommand{\calg}{\mathcal{G}}
\newcommand{\calh}{\mathcal{H}} 
 
\newcommand{\call}{\mathcal{L}} 
\newcommand{\caln}{\mathcal{N}}

\newcommand{\calx}{\mathcal{X}}

 \newcommand{\bfK}{{\mathbf K}}
\newcommand{\bfKH}{{\mathbf K}\mathbf{H}} \newcommand{\bfL}{{\mathbf L}}

\newcommand{\EGF}[2]{E_{#2}(#1)}               
\newcommand{\OrGF}[2]{\Or_{#2}(#1)}

\newcommand{\inj}{\operatorname{inj}}

\newcommand{\calcl}{{\calc\!\call}}

\newcommand{\Wh}{\operatorname{Wh}}

\newcommand{\chern}{\operatorname{chern}}

\newcommand{\HS}{\operatorname{HS}}

\newcommand{\Nil}{\operatorname{Nil}}

\newcommand{\edge}{\operatorname{edge}}

\newcommand{\class}{\operatorname{class}}

\newcommand{\KH}{K\!H}

\DeclareMathOperator{\cent}{cent}

\DeclareMathAlphabet{\matheurm}{U}{eur}{m}{n}

 \newcommand{\pt}{\text{pt}}

\newcommand{\Groupoids}{{\matheurm{Groupoids}}}
\newcommand{\Or}{\matheurm{Or}}

\renewcommand{\theenumi}{(\roman{enumi})}

\newcommand{\xycomsquare}[8]                   
{\xymatrix{#1 \ar[r]^{#2} \ar[d]^{#4} &
    #3 \ar[d]^{#5}  \\
    #6\ar[r]^{#7} & #8 }}

\begin{document}

\begin{abstract}
  We present the status of the Farrell-Jones Conjecture for algebraic
  $K$-theory for a group $G$ and arbitrary coefficient rings $R$.
  We add new groups for which the conjecture is known to be true and 
  study inheritance properties. We 
  discuss new applications, focussing on the Bass Conjecture, the
  Kaplansky Conjecture and conjectures generalizing Moody's Induction
  Theorem. Thus we extend the class of groups for which these
  conjectures are known considerably.
\end{abstract}



\maketitle

\setcounter{section}{-1}
\section{Introduction and statements of results}
\label{sec:introduction_and_statements_of_results}


\subsection{Background}
\label{subsec:Background_intro}

The \emph{Farrell-Jones Conjecture for algebraic $K$-theory} predicts
the structure of $K_n(RG)$ for a group $G$ and a ring $R$. There is
also an $L$-theory version. For applications in topology and geometry
the case $R = \IZ$ is the most important one since many topological
invariants of manifolds and $CW$-complexes such as the finiteness
obstruction, the Whitehead torsion and the surgery obstruction take
values in the algebraic $K$- or $L$-theory of the integral group ring
$\IZ \pi$ of the fundamental group $\pi$. The Farrell-Jones Conjecture
for $R = \IZ$ implies several famous conjectures, e.g., the
\emph{Novikov Conjecture}, (in high
dimensions) the \emph{Borel Conjecture},
and the triviality of compact $h$-cobordisms with torsionfree
fundamental group. On the other hand proofs of the
Farrell-Jones Conjecture for certain groups often rely on working with integral
coefficients since they are based on these geometric connections. This
is the reason why more is known about the algebraic $K$- and
$L$-theory of $\IZ G$ than of $\IC G$ which is in some sense
surprising since $\IC G$ has better ring theoretic properties than
$\IZ G$. For the status of the Farrell-Jones Conjecture with
coefficient in $\IZ$  we refer for instance to \cite[Sections~5.2
and~5.3]{Lueck-Reich(2005)}.

Recently the geometric approaches have been generalized so far that
they also apply to other coefficient rings than $\IZ$ (see for
instance Farrell-Linnell~\cite{Farrell-Linnell(2003b)},
Bartels-Reich~\cite{Bartels-Reich(2005JAMS)},
Bartels-L\"uck-Reich~\cite{Bartels-Lueck-Reich(2007hyper)},
Quinn~\cite{Quinn(2005)}).  This is interesting for algebraic and ring
theoretic applications, where one would like to consider for example fields, rings
of integers in algebraic number fields and  integral domains.
The purpose of this article is to describe the status of the
Farrell-Jones Conjecture for algebraic $K$-theory for arbitrary
coefficient rings and to discuss applications, for instance to the
\emph{Bass Conjectures}, the \emph{Kaplansky Conjecture},
generalizations of \emph{Moody's Induction Theorem}, \emph{Nil-groups} and
\emph{Fuglede-Kadison determinants}.


\subsection{Status of the Farrell-Jones Conjecture for algebraic $K$-theory}
\label{subsec:Status_of_the_Farrell-Jones_Conjecture_for_algebraic_K-theoryd_intro}

There is a stronger version of the Farrell-Jones Conjecture, the so
called \emph{Fibered Farrell-Jones Conjecture}. The Fibered
Farrell-Jones Conjecture does imply the Farrell-Jones Conjecture and
has better inheritance properties than the Farrell-Jones Conjecture.
We will give the precise technical formulations of these
conjectures in
Section~\ref{sec:The_Farrell-Jones_Conjecture_for_algebraic_K-Theory}.
The original source for the (Fibered) Farrell-Jones Conjecture is the paper by
Farrell-Jones~\cite[1.6 on page~257 and~1.7 on page~262]{Farrell-Jones(1993a)}.

Ring will always mean associative ring with unit. It is not
necessarily commutative. Fields are understood to be commutative
unless they are called skew-fields.

One of the main results of this article is the next theorem whose
proof will be given in Subsection~\ref{subsec:Extensions}.

\begin{theorem}
 \label{the:class_FJC}
 Let $R$ be a ring. Let $\calfj(R)$ be the class of groups which
 satisfy the Fibered Farrell-Jones Conjecture for algebraic $K$-theory
 with coefficients in $R$.  Then

\begin{enumerate}

\item \label{the:class_FJC:word_hyper_virt_nilpotent} Every
      word-hyperbolic group and every virtually nilpotent group
      belongs to $\calfj(R)$;

\item \label{the:class_FJC:products} If $G_1$ and $G_2$ belong to
      $\calfj(R)$, then $G_1 \times G_2$ belongs to $\calfj(R)$;

\item \label{the:class_FJC:colim} Let $\{G_i \mid i\in I\}$ be a
      directed system of groups (with not necessarily injective structure maps)
      such that $G_i \in \calf$ for $i \in I$.
      Then $\colim_{i \in I} G_i$ belongs to $\calfj(R)$;

\item \label{the:class_FJC:subgroups} If $H$ is a subgroup of $G$ and
      $G \in \calfj(R)$, then $H \in \calfj(R)$.
\end{enumerate}

\end{theorem}

If one restricts to lower and middle $K$-theory for torsionfree groups
and regular rings $R$, the Farrell-Jones Conjecture for algebraic
$K$-theory reduces to the following easier to understand conclusions
which are already very interesting in their own right.

Let $R$ be a ring and let $G$ be a group. Denote by $i \colon R \to RG$
the obvious inclusion.  Sending $(g , [P]) \in G \times K_0(R)$ to the
class of the $RG$-automorphism
$$
R[G] \otimes_R P \to R[G] \otimes_R P, \quad u \otimes x \mapsto
ug^{-1} \otimes x
$$
defines a map $\Phi \colon G/[G,G]\otimes_{\IZ} K_0(R) \to
K_1(RG)$.  Define the homomorphism
\begin{eqnarray}
A := \Phi \oplus K_1(i) \colon \left(G/[G,G] \otimes_{\IZ} K_0(R)\right) \oplus
K_1(R) \to K_1(RG).
\label{map_A}
\end{eqnarray}

Define $\Wh^R(G)$ to be the cokernel of $A$.
If $\widetilde{K}_0(R) = 0$ and the obvious map $R^{\times} \to K_1(R)$
is surjective, then $\Wh^R(G)$ coincides with $K_1(RG)/\langle (r\cdot
g) \mid r \in R^{\times}, g \in G\rangle$.  If $R = \IZ$, then
$\Wh^{\IZ}(G)$ is the classical \emph{Whitehead group} $\Wh(G)$ which
appears for instance in the $s$-cobordism theorem.

\begin{theorem} \label{thm:FJC_for_torsionfree_G_regular_R_and_middle_lower_K-theory}

  Let $R$ be a regular ring. Suppose that $G$ is torsionfree and that
  the Farrell-Jones Conjecture for algebraic $K$-theory with
  coefficients in $R$ holds for $G$. Then

\begin{enumerate}

\item $K_n(RG) = 0$ for $n \le -1$;

\item The change of rings map $K_0(R) \to K_0(RG)$ is bijective. In
      particular $\widetilde{K}_0(RG)$ is trivial if and only if
      $\widetilde{K}_0(R)$ is trivial;

\item The Whitehead group $\Wh^R(G)$ is trivial.

\end{enumerate}
\end{theorem}

The proof of
Theorem~\ref{thm:FJC_for_torsionfree_G_regular_R_and_middle_lower_K-theory}
can be found in~\cite[Conjecture~1.1 on page~652, Conjecture~1.1 on
page~657 and Corollary~2.3 on page~685]{Lueck-Reich(2005)}.

In particular we get for $R = \IZ$ that $\Wh(G) = 0$,
$\widetilde{K}_0(\IZ G) = 0$ and $K_n(\IZ G) = 0$ for $n \le -1$
holds, if the torsionfree group $G$ satisfies the Farrell-Jones
Conjecture with coefficient in $\IZ$.

These vanishing results have important geometric consequences. Namely, let $G$ be a
finitely presented group. Then $\Wh(G)$ vanishes if and only if every compact
$h$-cobordism of dimension $\ge 6$ with $G$ as fundamental group is trivial, and
$\widetilde{K}_0(\IZ G)$ vanishes if and only if every finitely dominated $CW$-complex
with $G$ as fundamental group is homotopy equivalent to a compact $CW$-complex.
The vanishing of $\widetilde{K}_0(\IZ G)$ implies that the group $G$ is
already of type FF (which is sometimes also called type FL in the literature)
if it is of type FP (see~\cite[Chapter VIII, Section~6]{Brown(1982)}).

The conclusions appearing in
Theorem~\ref{thm:FJC_for_torsionfree_G_regular_R_and_middle_lower_K-theory}
are known to be true for a torsionfree group $G$ provided that $G$
belongs to the class $\calfj(R)$ appearing in
Theorem~\ref{the:class_FJC}. Examples are torsionfree subgroups of a
finite product $\prod_{i=1}^r G_i$, where each group $G_i$ is
word-hyperbolic. Subgroups of products of word-hyperbolic groups are studied for instance
in~\cite{Bridson-Miller(2004)}.

More information about torsionfree groups $G$ for which
Theorem~\ref{thm:FJC_for_torsionfree_G_regular_R_and_middle_lower_K-theory}
is true in the case $R = \IZ$ can be found in~\cite[Theorem~5.5 on
page~722]{Lueck-Reich(2005)}. There results due to Aravinda, Farrell,
Hu, Jones, Linnell, and Roushon are listed.


\subsection{Homotopy $K$-theory and rings with finite characteristic}
\label{subsec:Homotopy_K-theory_and_rings_with_finite_characteristic_intro}

Sometimes one wants to consider special rings or is interested in
rational information only.  Then the Fibered Farrell-Jones Conjecture
is known to be true for more groups. This is illustrated by Theorem~\ref{the:class_FJ_KH_andFJ_N} 
below. 

For the definition and basic properties of homotopy
$K$-theory we refer to Weibel~\cite{Weibel(1989)}. For a positive
integer $N$ let $\IZ[1/N]$ be the subring of $\IQ$ consisting of
rational numbers $m/n$ with $m,n \in \IZ$ for which each
prime dividing $n$ divides $N$. For an abelian group $A$ let $A[1/N]$
be $A \otimes_{\IZ}\IZ[1/N]$.
A ring $R$ is said to be \emph{of finite characteristic} if there is
an integer $N \geq 2$ such that $N \cdot 1_R = 0$.
In this case the minimal positive integer with this property is called the \emph{characteristic of $R$}.

The proof of the next theorem will be given in
Subsection~\ref{subsec:Homotopy_K-theory_and_rings_with_finite_characteristic}.
The differences of the conclusions appearing in
Theorem~\ref{the:class_FJC} and the following theorem are that now
virtually nilpotent is replaced by elementary amenable in
assertion~\ref{the:class_FJ_KH_andFJ_N:elem_amen_word_hyper} and that
there are new assertions about extensions and actions on trees, namely
assertions~\ref{the:class_FJ_KH_andFJ_N:extensions}
and~\ref{the:class_FJ_KH_andFJ_N:actions_on_trees}.
The extension result is intriguing since elementary amenable groups and word-hyperbolic groups
form separate branches in Bridson's universe of groups
(see~\cite {Bridson(2006)}).

\begin{theorem} \label{the:class_FJ_KH_andFJ_N}
  Let $R$ be a ring.
  Consider the following assertions for a group
  $G$.

\begin{itemize}

\item [(\KH)]
      The group $G$ satisfies the Fibered Farrell-Jones
      Conjecture for homotopy $K$-theory with coefficients in $R$;

\item [(FC)]
      The ring $R$ has finite characteristic $N$.
      The Fibered Farrell-Jones Conjecture for algebraic $K$-theory
      for $G$ with coefficients in $R$ for both the families $\Fin$
      and $\VCyc$ is true after applying $- \otimes_{\IZ} \IZ[1/N]$ to
      the assembly map.

\end{itemize}

Let $\calfj_{\KH}(R)$ be the class of groups for which assertion (\KH) holds.
If $R$ has finite characteristic, then let
$\calfj_{FC}(R)$ be the class of groups for which
assertion (FC) is true. Let $\calf$ be
$\calfj_{FC}(R)$ or $\calfj_{\KH}(R)$. Then:

\begin{enumerate}

\item \label{the:class_FJ_KH_andFJ_N:elem_amen_word_hyper} Every
       word-hyperbolic and every elementary amenable group
      belongs to $\calf$;

\item \label{the:class_FJ_KH_andFJ_N:products} If $G_1$ and $G_2$
      belong to $\calf$, then $G_1 \times G_2$ belongs to $\calf$;

\item \label{the:class_FJ_KH_andFJ_N:colim} Let $\{G_i \mid i\in I\}$
      be a directed system of groups (with not necessarily injective structure maps)
      such that $G_i \in \calf$ for $i
      \in I$.  Then $\colim_{i \in I} G_i$ belongs to $\calf$;

\item \label{the:class_FJ_KH_andFJ_N:subgroups} If $H$ is a subgroup
      of $G$ and $G \in \calf$, then $H \in \calf$;

\item \label{the:class_FJ_KH_andFJ_N:extensions} Let $1 \to H \to G
      \to Q \to 1$ be an extension of groups such that $H$ is either
      elementary amenable or word-hyperbolic and $Q$ belongs to
      $\calf$. Then $G$ belongs to $\calf$;

\item \label{the:class_FJ_KH_andFJ_N:actions_on_trees}
      Suppose that $G$ acts on a tree $T$.
      Assume that for each $x \in T$ the isotropy
      group $G_x$ belongs to $\calf$. Then $G$ belongs to $\calf$.
\end{enumerate}

Moreover, if $R$ has finite characteristic then we have
$\calfj_{\KH}(R) \subseteq \calfj_{FC}(R)$.

\end{theorem}

\begin{corollary}
\label{cor:FJC_for_torsionfree_G_regular_R_and_middle_lower_K-theory_rationally}
Let $R$ be a regular ring of finite characteristic $N$.
Let $G$ be torsionfree.
Suppose that $G$
belongs to the class $\calfj_{FC}(R)$ defined in
Theorem~\ref{the:class_FJ_KH_andFJ_N}.
Then

\begin{enumerate}

\item $K_n(RG)[1/N] = 0$ for $n \le -1$;

\item The change of rings map induces a bijection $K_0(R)[1/N] \to
      K_0(RG)[1/N]$.  In particular $\widetilde{K}_0(RG)[1/N]$ is
      trivial if and only if $\widetilde{K}_0(R)[1/N]$ is trivial;

\item $\Wh^R(G)[1/N]$ is trivial.

\end{enumerate}
\end{corollary}

The proof of
Corollary~\ref{cor:FJC_for_torsionfree_G_regular_R_and_middle_lower_K-theory_rationally}
is analogous to the one of
Theorem~\ref{thm:FJC_for_torsionfree_G_regular_R_and_middle_lower_K-theory}.
Corollary~\ref{cor:FJC_for_torsionfree_G_regular_R_and_middle_lower_K-theory_rationally}
together with Theorem~\ref{the:class_FJ_KH_andFJ_N} substantially
extends Theorem~1.1 of Farrell-Linnell~\cite{Farrell-Linnell(2003b)},
where $\Wh^F(G) \otimes_{\IZ} \IQ = 0$ is proven for $G$ a torsionfree
elementary amenable group and $F$ a field of prime characteristic.

The vanishing of $\Wh^{\IF_p}(G)\otimes_{\IZ}\IQ$ is needed
in the definition of a $p$-adic logarithmic Fuglede-Kadison
determinant for $G$ due to Deninger~\cite{Deninger(2006)}.


\subsection{Induction from finite subgroups}
\label{subsec:Induction_from_finite_subgroupsd_intro}

The next result will be explained and proven in
Section~\ref{sec:The_projective_class_group_and_induction_from_finite_subgroups}.

\begin{theorem}
\label{the:FJC_implies_Moodys_induction}
Let $G$ be a group. Then
\begin{enumerate}
\item
\label{the:FJC_implies_Moodys_induction:reg_ring_Q_subset_R}
Let $R$ be a regular ring such that the order of any finite subgroup
of $G$ is invertible in $R$, e.g., a field of characteristic zero.
Suppose that $G$ satisfies the Farrell-Jones Conjecture with
coefficients in $R$. Then the map given by induction from finite
subgroups of $G$ (see \eqref{I(G,R):colim_FIN_to_K_0(RG)})
$$I(G,R) \colon \colim_{\OrGF{G}{\Fin}} K_0(RH) \to K_0(RG) $$
is an isomorphism;

\item
\label{the:FJC_implies_Moodys_induction:Field_char_p}
Let $D$ be a skew-field of characteristic $p$ for a prime number $p$.
Suppose that $G$ satisfies the Farrell-Jones Conjecture with
coefficients in $D$ after applying $- \otimes_{\IZ}\IZ[1/p]$.

Then the map
$$I(G,D)[1/p] \colon \colim_{\OrGF{G}{\Fin}} K_0(DH)[1/p] \to
K_0(DG)[1/p]$$
is an isomorphism.

\end{enumerate}
\end{theorem}

Theorem~\ref{the:FJC_implies_Moodys_induction} is an example of a result of the type that
certain $K$-groups of a group ring are given by induction over finite subgroups. A
prominent example is \emph{Moody's Induction Theorem} (see~\cite{Cliff-Weiss(1988)},
\cite{Moody(1987)}, \cite{Moody(1989)}) which deals with the surjectivity of a corresponding map
to $G_0(RG)$ instead of $K_0(RG)$ for virtually poly-$\IZ$ groups $G$ and
Noetherian rings $R$. If $R$ is regular and the order of any finite subgroup in the
virtually poly-$\IZ$ group $G$ is invertible in $R$, then $RG$ is regular and there is no
difference between the $G$-theoretic and the $K$-theoretic statement. Thus
Theorem~\ref{the:FJC_implies_Moodys_induction}~\ref{the:FJC_implies_Moodys_induction:reg_ring_Q_subset_R}
is linked to Moody's induction theorem.

Every group in the family $\calfj(R)$ appearing in Theorem~\ref{the:class_FJC} satisfies
the assumptions of
Theorem~\ref{the:FJC_implies_Moodys_induction}~\ref{the:FJC_implies_Moodys_induction:reg_ring_Q_subset_R}.

Theorem~\ref{the:FJC_implies_Moodys_induction}~\ref{the:FJC_implies_Moodys_induction:Field_char_p}
applies to every group $G$ which belongs to the family $\calfj_{FC}(D)$ appearing in
Theorem~\ref{the:class_FJ_KH_andFJ_N}. Thus we
have substantially generalized Theorem~1.2 of
Farrell-Linnell~\cite{Farrell-Linnell(2003b)}, where the surjectivity of the map
$$I(G,F)\otimes_{\IZ} \IQ \colon \colim_{\OrGF{G}{\Fin}} K_0(FH)\otimes_{\IZ} \IQ \to
K_0(FG)\otimes_{\IZ} \IQ$$
is proven for elementary amenable groups $G$ and fields $F$ of
prime characteristic.


\subsection{Bass Conjectures}
\label{subsec:Bass_Conjecturesd_intro}

The following conjecture is due to Bass~\cite[4.5]{Bass(1979)}.

\begin{conjecture}[Bass Conjecture for commutative integral domains]
\label{con:Bass_Conjecture_for_integral_domains}
Let $R$ be a commutative integral domain and let $G$ be a group. Let
$g \in G$. Suppose that either the order $|g|$ is
infinite or that the order $|g|$ is finite and not invertible in $R$.

Then for every finitely generated projective $RG$-module the value of
its Hattori-Stallings rank $\HS_{RG}(P)$ at $(g)$
(see~\eqref{HS_K_0(RG)_to_class_0}) is zero.
\end{conjecture}

The Bass Conjecture~\ref{con:Bass_Conjecture_for_integral_domains} can be interpreted
topologically. Namely, the Bass Conjecture~\ref{con:Bass_Conjecture_for_integral_domains}
is true for a finitely presented group $G$ in the case $R = \IZ$ if and only if every
homotopy idempotent selfmap of an oriented smooth closed manifold whose dimension is
greater than 2 and whose fundamental group is isomorphic to $G$ is homotopic to one that
has precisely one fixed point (see~\cite{Berrick-Chatterji-Mislin(2007)}). The Bass
Conjecture~\ref{con:Bass_Conjecture_for_integral_domains} for $G$ 
in the case $R = \IZ$ (or $R = \IC$) also implies for a
finitely dominated $CW$-complex with fundamental group $G$ that its Euler characteristic
agrees with the $L^2$-Euler characteristic of its universal covering
(see~\cite{Eckmann(1996b)}).

The next results follows from the argument in
\cite[Section~5]{Farrell-Linnell(2003b)}.

\begin{theorem} \label{the:FJC_implies_Bass_for_R_to_F_for_specific_values}
  Let $G$ be a group. Suppose that
  $$I(G,F)\otimes_{\IZ} \IQ \colon \colim_{\OrGF{G}{\Fin}} K_0(FH)
  \otimes_{\IZ} \IQ \to K_0(FG) \otimes_{\IZ} \IQ $$
  is surjective for
  all fields $F$ of prime characteristic.

  Then the Bass
  Conjecture~\ref{con:Bass_Conjecture_for_integral_domains}
  is satisfied for every commutative integral domain $R$.
\end{theorem}

Hence by Theorem~\ref{the:FJC_implies_Moodys_induction} the Bass
Conjecture for commutative integral
domains~\ref{con:Bass_Conjecture_for_integral_domains} is true for
every group $G$ which lies in the class $\calfj_{FC}(F)$ for all
fields $F$ of prime characteristic. The case of elementary amenable
groups has already been treated
in~\cite[Theorem~1.6]{Farrell-Linnell(2003b)}. The Bass
Conjecture~\ref{con:Bass_Conjecture_for_integral_domains} has been
proved by Bass~\cite[Proposition~6.2 and Theorem~6.3]{Bass(1979)} for $R =
\IC$ and $G$ a linear group,  by
Linnell~\cite[Lemma~4.1]{Linnell(1983b)} for $|g| < \infty$ and $R =
\IZ$, and by Eckmann~\cite[Theorem~3.3]{Eckmann(1986)} for $R = \IQ$
provided that $G$ has at most cohomological dimension $2$ over $\IQ$.
Further results are proved by Emmanouil~\cite{Emmanouil(1998a)}.

Here is another version of the Bass Conjecture.

\begin{conjecture}[Bass Conjecture for fields of characteristic zero as coefficients]
\label{con:Bass_Conjecture_for_FG}
Let $F$ be a field of characteristic zero and let $G$ be a group.  The
Hattori-Stallings homomorphism (see~\eqref{HS_K_0(RG)_to_class_0})
induces an isomorphism
$$\HS_{FG} \colon K_0(FG) \otimes_{\IZ} F \to \class_F(G)_f.$$
\end{conjecture}

For a field $F$ of characteristic zero Conjecture~\ref{con:Bass_Conjecture_for_FG}
obviously implies Conjecture~\ref{con:Bass_Conjecture_for_integral_domains}.
The proof of the next result will be given in
Section~\ref{sec:Bass_Conjectures}.

\begin{theorem}
\label{the:FJC_implies_BC_for_FG}
Let $F$ be a field of characteristic zero and let $G$ be a group. If
$G$ satisfies the Farrell-Jones Conjecture with coefficients in $F$,
then $G$ satisfies the
Bass-Conjecture~\ref{con:Bass_Conjecture_for_FG} for $F$.
\end{theorem}

In particular the Bass Conjecture for a field $F$ of characteristic
zero as coefficients~\ref{con:Bass_Conjecture_for_FG} is true for all
groups in the class $\calfj(F)$ by Theorem~\ref{the:class_FJC}.

Berrick-Chatterji-Mislin~\cite{Berrick-Chatterji-Mislin(2004)} prove
that a group $G$ satisfies the
Bass Conjecture~\ref{con:Bass_Conjecture_for_FG} for $F = \IC$ and the
Bass Conjecture~\ref{con:Bass_Conjecture_for_integral_domains} for integral
domains for $R = \IZ$ if $G$ satisfies the Bost Conjecture. Because
the Bost Conjecture is known for many groups, this is also true for
the Bass Conjecture for $F = \IC$. Since the Bost Conjecture deals
with $l^2$-spaces, this strategy can only work for subrings of $\IC$.

The following result describes another conclusion of the Farrell-Jones Conjecture, which is 
in the spirit of the
Bass Conjecture. It is hence true
for all groups in the family $\calfj(R)$.

\begin{theorem} \label{the:Bass_for_R_to_F}
  Let $G$ be a group. Let $R$ be a commutative integral domain with quotient field $F$
  such that no prime divisor of the order of a finite subgroup of $G$ is invertible in
  $R$. (An example is $R = \IZ$ and $F = \IQ$.)  Suppose that $G$ satisfies the
  Farrell-Jones Conjecture for algebraic $K$-theory with coefficients in $R$.

  Then the change of rings homomorphisms
  $$K_0(RG) \otimes_{\IZ} \IQ \to K_0(FG) \otimes_{\IZ} \IQ$$
  agrees
  with the composite
  $$K_0(RG) \otimes_{\IZ} \IQ \to K_0(R) \otimes_{\IZ} \IQ \to K_0(F)
  \otimes_{\IZ} \IQ \to K_0(FG) \otimes_{\IZ} \IQ$$
  where the three
  maps come from the change of ring homomorphisms associated to the
  augmentation $RG \to R$, the inclusion $R \to F$ and the inclusion
  $F \to FG$. In particular the homomorphism
  $$\widetilde{K}_0(RG)\otimes_{\IZ} \IQ \to \widetilde{K}_0(FG)\otimes_{\IZ} \IQ$$
  is trivial.
\end{theorem}

(We remind the reader, that $\widetilde{K}_0(RG)$ and
$\widetilde{K}_0(FG)$ respectively is the cokernel of
the canonical map $K_0(\IZ) \to K_0(RG)$ and the canonical map $K_0(\IZ) \to K_0(FG)$
respectively.)

If $G$ is finite, $G$ satisfies the Farrell-Jones Conjecture for
algebraic $K$-theory with coefficients in $R$ for trivial reasons and
hence Theorem~\ref{the:Bass_for_R_to_F} reduces to a Theorem of Swan
(see~\cite[Theorem~8.1]{Swan(1960a)},~\cite[Corollary~4.2]{Bass(1979)}).

The conclusion of Theorem~\ref{the:Bass_for_R_to_F} is related to the
theorem that for every group $G$ the change of rings homomorphism $\widetilde{K}_0(\IZ G) \to
\widetilde{K}_0(\caln(G))$
is trivial, where $\caln(G)$ is the group von Neumann algebra
(see~\cite[Theorem~9.62 on page~363]{Lueck(2002)},~\cite {Schafer(2000)}).


\subsection{The Kaplansky Conjecture}
\label{subsec:The_Kaplansky_Conjecture_intro}

\begin{conjecture}[Kaplansky Conjecture]\label{con:Kaplansky_Conjecture}
  Let $R$ be an integral domain and let $G$ be a torsionfree group.
  Then all idempotents of $RG$ are trivial, i.e., equal to $0$ or $1$.
\end{conjecture}

In the next theorem we will use the notion of a \emph{sofic group} that was introduced by
Gromov and originally called \emph{subamenable group}.  Every residually amenable group is
sofic but the converse is not true.  The class of sofic groups is closed under taking
subgroups, direct products, free amalgamated products, colimits and inverse limits, and,
if $H$ is a sofic normal subgroup of $G$ with amenable quotient $G/H$, then $G$ is sofic.
For more information about the notion of a sofic group we refer
to~\cite{Elek-Szabo(2006)}.

The proof and further explanations of the next theorem will be given in
Section~\ref{sec:The_Kaplansky_Conjecture}.

\begin{theorem} \label{the:FJC_and_Kaplansky}
  Let $G$ be a group. Let $R$ be a ring whose idempotents are all trivial. Suppose that
  $$K_0(R) \otimes_{\IZ} \IQ  \xrightarrow{} K_0(RG) \otimes_{\IZ} \IQ$$
  is an isomorphism.

  Then the Kaplansky Conjecture holds for $R$ and $G$ if one of the
  following conditions is satisfied:

\begin{enumerate}

\item $RG$ is stably finite;

\item $R$ is a field of characteristic zero;

\item $R$ is a skew-field and $G$ is sofic.

\end{enumerate}

\end{theorem}

Next we discuss some special cases of
Theorem~\ref{the:FJC_and_Kaplansky}.
Notice that we get assertions also for skew-fields and not only for fields
as coefficients.

Theorem~\ref{thm:FJC_for_torsionfree_G_regular_R_and_middle_lower_K-theory}
and Theorem~\ref{the:FJC_and_Kaplansky} imply that for a skew-field $D$ of
characteristic zero  and a torsionfree  group $G$ belonging to the class
of groups $\calfj(D)$ defined in Theorem~\ref{the:class_FJC}
the Kaplansky Conjecture~\ref{con:Kaplansky_Conjecture} is true
for $DG$, provided that $D$ is commutative or that $G$ is sofic.

Suppose that $D$ is a skew-field of prime characteristic $p$, all
finite subgroups of $G$ are $p$-groups and $G$ belongs to the class
$\calfj_{FC}(D)$ defined in Theorem~\ref{the:class_FJ_KH_andFJ_N}.
Then $K_0(D) \otimes_{\IZ} \IQ \xrightarrow{\cong}
K_0(DG) \otimes_{\IZ} \IQ$ is an isomorphism by
Theorem~\ref{the:FJC_implies_Moodys_induction}~\ref{the:FJC_implies_Moodys_induction:Field_char_p}
since for a finite $p$-group $H$ the group ring $DH$ is a local ring
and hence $\widetilde{K}_0(DH) = 0$.  If we furthermore assume that
$G$ is sofic, then Theorem~\ref{the:FJC_and_Kaplansky} implies that all
idempotents in $DG$ are trivial.  This has already been proved in the
case where $G$ is elementary amenable and $D$ is commutative
by Farrell-Linnell~\cite[Theorem~1.7]{Farrell-Linnell(2003b)}.

To the authors' knowledge there is no example of a group which is not sofic and of a
group which is word-hyperbolic and not residually finite in the literature.
So it is conceivable that all word-hyperbolic groups are sofic.

Next we mention some results of others.

Let $F$ be a field of characteristic zero and $u = \sum_{g \in G} x_g
\cdot g\in K[G]$ be an idempotent. Let $K$ be the finitely generated
field extension of $\IQ$ given by $K = \IQ(x_g \mid g \in G)$. Obviously
$u \in KG$ is an idempotent. There exists an embedding of $K$ in $\IC$. Hence all
idempotents in $FG$ are non-trivial if all idempotents in $\IC G$ are
trivial.

The \emph{Kadison Conjecture} says that all idempotents in the
reduced group $C^*$-algebra $C^*_r(G)$ of a torsionfree group are
trivial.  Hence the Kadison Conjecture implies the Kaplansky
Conjecture for all fields of characteristic zero. The Kadison
Conjecture follows from the \emph{Baum-Connes Conjecture} (as
explained for instance in~\cite{Lueck(2002d)}, \cite[1.8.1
and~1.8.2]{Lueck-Reich(2005)}). Hence a torsionfree group $G$
satisfies the Kaplansky Conjecture for all fields of
characteristic zero, if it satisfies the Baum-Connes Conjecture.
For a survey of groups satisfying the Baum-Connes Conjecture we
refer to~\cite[Sections~5.1 and~5.3]{Lueck-Reich(2005)}. We
mention that subgroups of word-hyperbolic groups satisfy the
Baum-Connes Conjecture by a result of
Mineyev-Yu~\cite[Theorem~20]{Mineyev-Yu(2002)} based on the work
of Lafforgue~\cite{Lafforgue(1998)} (see
also~\cite{Skandalis(1999)}).  A proof of the Kadison Conjecture
for a torsionfree word-hyperbolic group using cyclic homology is
given by Puschnigg~\cite{Puschnigg(2002)}. Notice that all these
analytic methods do only work for fields of characteristic zero
and cannot be extended to skew-fields or fields of prime
characteristic.

Formanek~\cite{Formanek(1973)} (see also~\cite[Lemma~4.1 and Proposition~4.2]{Burger-Valette(1998)})
has shown that all idempotents of $FG$ are trivial
provided that $F$ is a field of prime characteristic $p$, the group
$G$ contains no $p$-torsion and there do not exist an element $g \in
G, g \not= 1$ and an integer $k \ge 1$ such that $g$ and $g^{p^k}$ are
conjugate, or,  provided that $F$ is a field of characteristic zero and
there are infinitely many primes $p$ for which there do not exist an
element $g \in G, g \not= 1$ and an integer $k \ge 1$ such that $g$
and $g^{p^k}$ are conjugate.  Torsionfree word-hyperbolic groups
satisfy these conditions.  Hence Formanek's results imply that all
idempotents in $FG$ are trivial if $G$ is torsionfree word-hyperbolic
and $F$ is a field.

Delzant~\cite{Delzant(1997)} has proven the Kaplansky
Conjecture~\ref{con:Kaplansky_Conjecture} for all integral domains
$R$ for a torsionfree word-hyperbolic group $G$ provided that $G$
admits an appropriate action with large enough injectivity radius.
Delzant actually deals with zero-divisors and units as well.


\subsection{Homotopy invariance of $L^2$-torsion}
\label{subsec:Homotopy_Invariance_of_L2-torsion_intro}

The following conjecture for a group $G$ is stated and explained in
L\"uck~\cite[Conjecture~3.94 (1) on page~163]{Lueck(2002)}.

\begin{conjecture}
\label{con:hom_inv_o_L2-torsion}
Define the homomorphism
$$\Phi = \Phi^G\colon \Wh(G) \to \IR$$
by sending the class $[A]$ of an invertible matrix $A \in GL_n(\IZ
G)$ to $\ln(\det(r_A^{(2)}))$, where $\det(r_A^{(2)})$ is the
Fuglede-Kadison determinant of the $G$-equivariant bounded
operator $l^2(G)^n \to l^2(G)^n$ given by right multiplication
with $A$.

Then $\Phi$ is trivial.
\end{conjecture}

It is important because of the following conclusion explained in
\cite[Theorem~3.94~(1)) on page~161]{Lueck(2002)}: If $X$ and $Y$ are
det-$L^2$-acyclic finite $G$-$CW$-complexes, which are $G$-homotopy
equivalent, then their \emph{$L^2$-torsions} agree:
$$\rho^{(2)}(X;\caln(G)) = \rho^{(2)}(Y;\caln(G)).$$

\begin{theorem}\label{the:FJC_homotopy_inv_of_L2_torsion}
  Suppose that $G$ satisfies the Farrell-Jones Conjecture for algebraic
  $K$-theory with coefficients in $\IZ$. Then $G$ satisfies
  Conjecture~\ref{con:hom_inv_o_L2-torsion}.
\end{theorem}

We will omit the proof of
Theorem~\ref{the:FJC_homotopy_inv_of_L2_torsion} since it is similar
to the one of Theorem~\ref{the:Bass_for_R_to_F} using the fact that
for a finite group $H$ we have $\widetilde{K}_0(\IZ H) \otimes_{\IZ} \IQ = 0$
and Conjecture~\ref{con:hom_inv_o_L2-torsion} is true
for finite groups for elementary reasons.

Let $G$ be a torsionfree word-hyperbolic group. Suppose that its
$L^2$-Betti numbers $b_p^{(2)}(G)$ are trivial for all $p \ge 0$ and
that $G$ is of $\det \ge 1 $-class. (If $G$ is residually finite, it
is of $\det \ge 1$- class.) Choose a cocompact model for $EG$. Then we
can define the \emph{$L^2$-torsion of $G$}
$$\rho^{(2)}(G) := \rho^{(2)}(EG;\caln(G)) \in \IR.$$
This is indeed
independent of the choice of a cocompact model for $EG$ and hence depends
only on $G$ by Theorem~\ref{the:class_FJC} and
Theorem~\ref{the:FJC_homotopy_inv_of_L2_torsion}.  If $M$ is a closed
hyperbolic manifold of dimension $2n+1$, then its fundamental group
$\pi = \pi_1(M)$ satisfies all these assumptions and there exists a
number $C_n > 0$ depending only on $n$ such that $C_n \cdot (-1)^n
\cdot \rho^{(2)}(\pi)$ is the volume of $M$.
Hence $\rho^{(2)}(G)$ can
be viewed as a kind of volume of a word-hyperbolic group $G$,
provided $G$ satisfies the above assumptions.


\subsection{Searching for counterexamples}
\label{subsec:Searching_for_conterexamplesd_intro}
There is no group known for which the Farrell-Jones Conjecture, the
Fibered Farrell-Jones Conjecture or the Baum-Connes Conjecture is
false. However, Higson, Lafforgue and
Skandalis~\cite[Section~7]{Higson-Lafforgue-Skandalis(2002)} construct
counterexamples to the \emph{Baum-Connes-Conjecture with
coefficients}, actually with a commutative $C^*$-algebra as
coefficients. They describe precisely what properties a group $\Gamma$
must have so that it does \emph{not} satisfy the Baum-Connes
Conjecture with coefficients.  Gromov~\cite{Gromov(2000)} describes
the construction of such a group $\Gamma$ as a colimit over a directed
system of groups $\{G_i \mid i \in I\}$ such that each $G_i$ is
word-hyperbolic.  We conclude from Bartels-Echterhoff-L\"uck~\cite{Bartels-Echterhoff-Lueck(2007colim)}
and Theorem~\ref{the:class_FJC}
that the Fibered Farrell-Jones Conjecture and the Bost Conjecture do hold for
$\Gamma$.


\subsection{Nil-groups}
\label{subsec:Nil-groupsd_intro}
In Section~\ref{sec:Nil-groups} we discuss some consequences for
Nil-groups in the sense of Bass and Waldhausen and for the passage from algebraic $K$-theory
to homotopy $K$-theory.
There we will prove the following  application of
Theorem~\ref{the:class_FJC}~\ref{the:class_FJC:word_hyper_virt_nilpotent}
to Waldhausen's Nil-groups.

\begin{theorem}
\label{thm:waldhausen-nil-is-rationally-trivial}
Let $G$, $H$ and $K$ be finite groups. Let $C := \IZ[K]$.
\begin{enumerate}
\item \label{thm:waldhausen-nil:amalgamated}
      Let $\alpha \colon K \to G$ and $\beta \colon K \to H$ be injective
      group homeomorphisms.
      Consider the $C$ bimodules $A' := \IZ[G - \alpha(K)]$ and $B' := \IZ[H - \beta(K)]$.
      Then
      \[
      \Nil_n(C; A',B') \otimes_\IZ \IQ = 0;
      \]
\item \label{thm:waldhausen-nil:laurent}
      Let $\alpha \colon K \to G$ and $\beta \colon K \to G$ be injective
      group homeomorphisms.
      Let $A' := \IZ[G - \alpha(K)]$, $A'' := \IZ[H - K]$ and $A := \IZ[G]$.
      Then
      \[
      \Nil_n(C;{_{\alpha}A'_{\alpha}},{_{\beta}A''_{\beta}},{_{\beta}A_{\alpha}},{_{\alpha}A_{\beta}})
      \otimes_\IZ \IQ = 0,
      \]
      where we used the lower indices to indicate the relevant $C$-bimodule structures.
\end{enumerate}
\end{theorem}


\subsection{The Farrell-Jones Conjecture for L-theory}
\label{subsec:The_Farrell-Jones_Conjecture_for_L-theoryd_intro}
In Section~\ref{sec:The_Farrell-Jones_Conjecture_for_L-theory} we
briefly explain some results about the $L$-theoretic version of the
Farrell-Jones Conjecture.


\section{Inheritance Properties of the (Fibered) Isomorphism Conjecture}
\label{sec:Inheritance_Properties_of_the_(Fibered)_Isomorphism_Conjecture}

In this section we formulate a (Fibered) Isomorphism Conjecture
for a given equivariant homology theory and a family of groups.
In this general setting we study the behavior of this Fibered
Isomorphism Conjecture under directed colimits and extensions. The
Farrell-Jones Conjecture is a special case, one has to choose a
specific equivariant homology theory and a specific family of subgroups.  The
payoff of this general setting is that some of the proofs become
easier and more transparent and that it applies to other related
Isomorphism Conjectures such as the \emph{Farrell-Jones Conjecture
for $L$-theory}, the \emph{Baum-Connes
  Conjecture}, the \emph{Bost Conjecture} and other types of
Isomorphism Conjectures predicting the bijectivity of certain assembly
maps.

For this section we fix the following data

\begin{itemize}

\item a discrete group $G$;

\item an equivariant homology theory $\calh^?_*$ with values in
      $\Lambda$-modules;

\item a class of groups $\calc$ closed under isomorphisms, taking
      subgroups and taking quotients, e.g., the family $\Fin$ of
      finite groups and the family $\VCyc$ of virtually cyclic groups.
      For a group $G$ we denote by $\calc(G)$ the family of subgroups of
      $G$ which belong to $\calc$.

\end{itemize}

Here \emph{equivariant homology theory} with values in
$\Lambda$-modules for a commutative ring $\Lambda$ satisfying the
disjoint union axiom is understood in the sense of
\cite[Section~1]{Lueck(2002b)} with one important modification: We
require that for every group homomorphism $\alpha \colon H \to K$ we
get a natural transformation
$$\ind_{\alpha} \colon \calh_*^H(X) \to \\\calh^K_n(\ind_{\alpha} X)$$
satisfying the obvious variations of the axioms a.) Compatibility with
the boundary operator b.)  Naturality and c.) Compatibility with
conjugation, but the map $\ind_{\alpha}$ is only required to be an
isomorphism in the case, where $X = \pt$ and $\alpha$ is injective.
This implies that $\ind_{\alpha}$ is bijective for a $G$-$CW$-complex $X$
if the kernel of $\alpha$ acts freely on $X$
(see~\cite[Lemma~1.5]{Bartels-Echterhoff-Lueck(2007colim)}).
Every $\Groupoids$-spectrum gives an equivariant homology theory with values
in $\IZ$-modules in the sense above
(see~\cite[6.5]{Lueck-Reich(2005)}).  In particular we get an
equivariant homology theory with values in $\IZ$-modules in the sense
above for algebraic $K$-theory (see~\cite[Section 2]{Davis-Lueck(1998)},
\cite[Theorem~6.1]{Lueck-Reich(2005)}). If $\calh^?_*$ is an
equivariant homology theory with values in $\IZ$-modules, then
$\calh^?_* \otimes_{\IZ} \Lambda$ is an equivariant homology theory with
values in $\Lambda$-modules for $\IZ \subseteq \Lambda \subseteq \IQ$.

Notice that this is one of the key differences between the
Farrell-Jones Conjecture for algebraic $K$- and $L$-theory and the
Baum-Connes Conjecture for  topological $K$-theory of reduced group
$C^*$-algebras. In the latter case induction is only defined if the
kernel of the group homomorphism acts freely because the corresponding
spectrum lives over $\Groupoids^{\inj}$ and not over $\Groupoids$ as
in the Farrell-Jones setting (see~\cite[6.5]{Lueck-Reich(2005)}).


\subsection{The Fibered Isomorphism Conjecture for equivariant homology theories}
\label{subsec:The_Fibered_Isomorphism_Conjecture_for_equivariant_homology_theories}

A \emph{family of subgroups of $G$} is a collection of subgroups of
$G$ which is closed under conjugation and taking subgroups.  Let
$\EGF{G}{\calf}$ be the \emph{classifying space associated to
  $\calf$}.  It is uniquely characterized up to $G$-homotopy by the
properties that it is a $G$-$CW$-complex and that $\EGF{G}{\calf}^H$
is contractible if $H \in \calf$ and is empty if $H \notin \calf$.
For more information about these spaces $\EGF{G}{\calf}$ we refer to
the survey article~\cite{Lueck(2005s)}.  Given a group homomorphism
$\phi \colon K \to G$ and a family $\calf$ of subgroups of $G$, define
the family $\phi^*\calf$ of subgroups of $K$ by
\begin{eqnarray}
\phi^*\calf  & = & \{H \subseteq K \mid \phi(H) \in \calf\}.
\label{phiastcalf}
\end{eqnarray}
If $\phi$ is an inclusion of a subgroup, we also write $\calf|_K$
instead of $\phi^*\calf$.

\begin{definition}[(Fibered) Isomorphism Conjecture for $\calh^?_*$]
\label{def:(Fibered)_Isomorphism_Conjectures_for_calh?_ast}
A group $G$ together with a family of subgroups $\calf$ satisfies the
\emph{Isomorphism Conjecture for $\calh^?_*$} if the projection $\pr \colon
\EGF{G}{\calf} \to \pt$ to the one-point-space $\pt$ induces an
isomorphism
$$\calh^G_n(\pr) \colon \calh^G_n(\EGF{G}{\calf}) \xrightarrow{\cong}
\calh^G_n(\pt)$$
for $n \in \IZ$.

The pair $(G,\calf)$ satisfies the \emph{Fibered Isomorphism
  Conjecture for $\calh^?_*$} if for every group homomorphism $\phi \colon K \to G$ the
pair $(K,\phi^*\calf)$ satisfies the Isomorphism Conjecture.
\end{definition}

We mostly work with a fixed equivariant homology theory $\calh^?_{\ast}$
and hence we will often omit it in the statements.

The following results are proven in~\cite[Lemma~1.6]{Bartels-Lueck(2004ind)}
and~\cite[Lemma~1.2 and Theorem~2.4]{Bartels-Lueck(2006)}

\begin{lemma} 
\label{lem:enlarging_the_family_for_the_Fibered_Isomorphism_Conjecture}
  Let $G$ be a group and let $\calf \subset \calg$ be families of
  subgroups of $G$.  Suppose that $(G,\calf)$ satisfies the Fibered
  Isomorphism Conjecture.

  Then $(G,\calg)$ satisfies the Fibered Isomorphism Conjecture.
\end{lemma}

\begin{lemma} 
\label{lem:basic_inheritance_property_of_fibered_conjecture}
  Let $\phi \colon K \to G$ be a group homomorphism and let $\calf$ be
  a family of subgroups.  If $(G,\calf)$ satisfies the Fibered
  Isomorphism Conjecture, then $(K,\phi^*\calf)$ satisfies the Fibered
  Isomorphism Conjecture.
\end{lemma}

\begin{theorem}[Transitivity Principle] 
\label{the:transitivity}
  Let $\calf \subseteq \calg$ be two families of subgroups of $G$.
  Assume that for every element $H \in \calg$ the group $H$ satisfies
  the (Fibered) Isomorphism Conjecture for $\calf|_H$.

  Then $(G,\calg)$ satisfies the (Fibered) Isomorphism Conjecture if
  and only if $(G,\calf)$ satisfies the (Fibered) Isomorphism
  Conjecture.
\end{theorem}

The next lemma follows from  Lemma~\ref{lem:basic_inheritance_property_of_fibered_conjecture}
applied to the inclusion $H \to G$ since $\calc(H) = \calc(G)|_H$.

\begin{lemma} \label{lem:Fibered_Conjecture_passes_to_subgroups}
  Suppose that the Fibered Isomorphism Conjecture holds for $(G,\calc(G))$.
  Let $H \subseteq G$ be a subgroup.

  Then the Fibered Isomorphism Conjecture holds for $(H,\calc(H))$.
\end{lemma}


\subsection{Colimits over directed systems of groups}
\label{subsec:Colimits_over_directed_systems_of_groups_for_calh}

We collect some basic facts about the behavior of the Fibered
Isomorphism Conjecture under directed colimits.

We consider a \emph{directed} set $I$ and a \emph{directed system of
  groups} $\{G_i \mid i \in I\}$. The \emph{structure maps}
$\phi_{i,j} \colon G_i \to G_j$ for $i,j \in I$ with $i \le j$ are
\emph{not} required to be injective.  Let $\colim_{i \in I} G_i$ be
the \emph{colimit}.  Denote by $\psi_i \colon G_i \to G$ the
\emph{structure maps of the colimit} for $i \in I$.

We say that $G$ is the \emph{directed union of the subgroups} $\{G_i
\mid i \in I\}$ if $I$ is a directed set and $\{G_i \mid i \in I\}$ is
a directed system of subgroups, directed by inclusion, such that $G =
\bigcup_{i \in I} G_i$. This is essentially the same as a directed
system of groups such that all structure maps $\phi_{i,j}$ are
inclusions of groups and $G = \colim_{i \in I} G_i$.
For a group homomorphism $\psi \colon G' \to G$ define the $\Lambda$-map 
\[
\alpha_n(\psi) \colon \calh_n^{G'}(\pt) \to  \calh_n^G(\pt)
\]
as the composition of $\ind_{\psi}$ with the map
induced by the projection $\psi^* \pt \to \pt$
of $G$-spaces. 

The next definition is an extension of
\cite[Definition~3.1]{Bartels-Lueck(2004ind)}.

\begin{definition}[(Strongly) Continuous equivariant homology theory]
\label{def:continuous_equivariant_homology_theory} An equivariant
homology theory $\calh^?_*$ is called \emph{continuous} if for each
group $G$ which is the directed union of subgroups $\{G_i \mid i \in
I\}$ the $\Lambda$-map
$$\colim_{i \in I} \alpha_n(G_i \to G) \colon \colim_{i \in I}
\calh^{G_i}_n(\pt) \to \calh^G_n(\pt)$$
is an isomorphism for every $n
\in \IZ$.

An equivariant homology theory $\calh^?_*$ is called \emph{strongly
  continuous} if for each directed system of groups $\{G_i \mid i \in
I\}$ with $G = \colim_{i \in I} G_i$ the $\Lambda$-map
$$\colim_{i \in I} \alpha_n(\psi_i)\colon \colim_{i \in I}
\calh^{G_i}_n(\pt) \to \calh^G_n(\pt)$$
is an isomorphism for every $n
\in \IZ$.
\end{definition}

The next theorem generalizes the result of
Farrell-Linnell~\cite[Theorem~7.1]{Farrell-Linnell(2003a)} to a more
general setting about equivariant homology theories as developed in
Bartels-L\"uck~\cite{Bartels-Lueck(2004ind)}. Its proof can be found in
\cite[Theorems~3.4 and 4.6]{Bartels-Echterhoff-Lueck(2007colim)}.

\begin{theorem}
\label{the:isomorphism_conjecture_is_stable_under_colim}

\begin{enumerate}
\item
\label{the:isomorphism_conjecture_is_stable_under_colim:injective}
Let $G$ be the directed union $G = \bigcup_{i \in I} G_i$ of subgroups
$G_i$ Suppose that $\calh^?_*$ is continuous and that the (Fibered)
Isomorphism Conjecture is true for $(G_i,\calc(G_i))$ for all $i \in
I$.

Then the (Fibered) Isomorphism Conjecture is true for $(G,\calc(G))$;

\item
\label{the:isomorphism_conjecture_is_stable_under_colim:general}
Let $\{G_i \mid i \in I\}$ be a directed system of groups. Put $G =
\colim_{i \in I} G_i$.  Suppose that $\calh^?_*$ is strongly
continuous and that the Fibered Isomorphism Conjecture is true for
$(G_i,\calc(G_i))$ for all $i \in I$.

Then the Fibered Isomorphism Conjecture is true for $(G,\calc(G))$.
\end{enumerate}
\end{theorem}


\subsection{Extensions}
\label{subsec:Extensions_for_calh}

For the remainder of this section fix the following data:

\begin{itemize}

\item a discrete group $G$;

\item an equivariant homology theory $\calh^?_*$ with values in
      $\Lambda$-modules;

\item a class of groups $\calc$ closed under isomorphisms, taking
      subgroups and taking quotients, e.g., $\Fin$ or  $\VCyc$;

\item an exact sequence of groups $1 \to K \xrightarrow{i} G
      \xrightarrow{p} Q \to 1$.

\end{itemize}

We want to investigate the inheritance properties of the (Fibered)
Isomorphism
Conjecture~\ref{def:(Fibered)_Isomorphism_Conjectures_for_calh?_ast}
under exact sequences.

\begin{lemma} \label{lem:extensions_unfibered_conclusion}
  Suppose that the Fibered Isomorphism Conjecture holds for
  $(Q,\calc(Q))$ and for every $H \in p^*\calc(Q)$ the Isomorphism
  Conjecture is true for $(H,\calc(H))$.

  Then the Isomorphism Conjecture is true for $(G,\calc(G))$.
\end{lemma}
\begin{proof}
  Since  the Fibered Isomorphism Conjecture holds for
  $(Q,\calc(Q))$ by assumption, the Isomorphism Conjecture holds for
  $(G,p^*\calc(Q))$.  We have to show that the Isomorphism Conjecture
  holds for $(G,\calc(G))$. But this follows from the Transitivity
  Principle~\ref{the:transitivity}.
\end{proof}

\begin{lemma} \label{lem:extensions_fibered_conclusion}
  Suppose that the Fibered Isomorphism Conjecture holds for
  $(Q,\calc(Q))$. Then the following assertions are equivalent.

\begin{enumerate}

\item
      \label{lem:extensions_fibered_conclusion:FIC_for_calc_extensions_of_K}
      The Fibered Isomorphism Conjecture is true for
      $(p^{-1}(H),\calc(p^{-1}(H)))$ for every $H \in \calc(Q)$;

\item \label{lem:extensions_fibered_conclusion:FIC_for_G}

      The Fibered Isomorphism Conjecture is true for $(G,\calc(G))$.

\end{enumerate}
\end{lemma}
\begin{proof}\ref{lem:extensions_fibered_conclusion:FIC_for_G}
  $\Rightarrow $\ref{lem:extensions_fibered_conclusion:FIC_for_calc_extensions_of_K}
  This follows from
  Lemma~\ref{lem:basic_inheritance_property_of_fibered_conjecture}
  applied to the inclusion $p^{-1}(H) \to G$.
  \\[1mm]\ref{lem:extensions_fibered_conclusion:FIC_for_calc_extensions_of_K}
  $\Rightarrow $\ref{lem:extensions_fibered_conclusion:FIC_for_G} Let
  $q \colon L \to G$ be a group homomorphism. We have to show that
  $(L,q^*\calc(G))$ satisfies the Isomorphism Conjecture.  Since
  $(Q,\calc(Q))$ satisfies the Fibered Isomorphism Conjecture, we
  conclude that $(L,q^*p^*\calc(Q))$ satisfies the Isomorphism
  Conjecture. Because of the Transitivity principle~\ref{the:transitivity} it remains to
  show for any $H \subseteq L$ for which there exists $V \in \calc(Q)$
  with $q(H) \subseteq p^{-1}(V)$ that $(H,(q^*\calc(G))|_H)$
  satisfies the Isomorphism Conjecture. This follows from the
  assumption that $(p^{-1}(V), \calc(p^{-1}(V)))$ satisfies the
  Fibered Isomorphism Conjecture since the families $(q^*\calc(G))|_H$
  and $(q|_H)^*\calc(p^{-1}(V))$ coincide.
\end{proof}

\begin{lemma} \label{lem:extensions_with_finite_kernel}
  Suppose that $p^{-1}(H)$ belongs to $\calc(G)$ if $H \in \calc(Q)$.
  Then $(G,\calc(G))$ satisfies the Fibered Isomorphism Conjecture, if
  $(Q,\calc(Q))$ satisfies the Fibered Isomorphism Conjecture.
\end{lemma}
\begin{proof}
  This follows from
  Lemma~\ref{lem:basic_inheritance_property_of_fibered_conjecture}
  since $p^*\calc(Q) = \calc(G)$.
\end{proof}
Lemma~\ref{lem:extensions_with_finite_kernel} is interesting in the
case, where $\calc$ is $\Fin$ or $\VCyc$ and $K$ is finite.

\begin{lemma} \label{lem:products}

\begin{enumerate}

\item \label{lem:products:closed_under_products} Suppose that $H_1
      \times H_2$ belongs to $\calc$ if $H_1$ and $H_2$ belong to
      $\calc$. Then $(G_1 \times G_2,\calc(G_1 \times G_2))$ satisfies
      the Fibered Isomorphism Conjecture if and only if both
      $(G_1,\calc(G_1))$ and $(G_2,\calc(G_2))$ satisfy the Fibered
      Isomorphism Conjecture;

\item \label{lem:products:VCcy} Suppose that $(D_{\infty}\ \times
      D_{\infty},\VCyc(D_{\infty}\ \times D_{\infty}))$ satisfies the
      Fibered Isomorphism Conjecture, where $D_{\infty} = \IZ
      \rtimes\IZ/2$ is the infinite dihedral group.

      Then $(G_1 \times G_2,\VCyc(G_1 \times G_2))$ satisfies the
      Fibered Isomorphism Conjecture if and only if both
      $(G_1,\VCyc(G_1))$ and $(G_1,\VCyc(G_1))$ satisfy the Fibered
      Isomorphism Conjecture.

\end{enumerate}
\end{lemma}
\begin{proof}\ref{lem:products:closed_under_products} If the Fibered Isomorphism
  Conjecture holds for $(G_1 \times G_2,\calc(G_1) \times
  \calc(G_2))$, it holds by
  Lemma~\ref{lem:Fibered_Conjecture_passes_to_subgroups} also for
  $(G_i,\calc(G_i))$ since $\calc(G_1 \times G_2)|_{G_i} =
  \calc(G_i)$.

  Suppose that the Fibered Isomorphism Conjecture holds for both
  $(G_1,\calc(G_1))$ and $(G_2,\calc(G_2))$.  In view of
  Lemma~\ref{lem:extensions_fibered_conclusion} we can assume without
  loss of generality that $G_2$ belongs to $\calc$. Applying this
  argument again, we can assume without loss of generality that $G_1$
  and $G_2$ belong $\calc$.  This case is obviously true since $G_1
  \times G_2 \in \calc$.
  (Compare \cite[Lemma~5.1]{Roushon(2006)} for a similar argument.)
  \\[1mm]\ref{lem:products:VCcy} Analogously to the proof of
  assertion~\ref{lem:products:closed_under_products}, one reduces the
  claim to the assertion that $(G_1 \times G_2, \VCyc(G_1 \times G_2))$
  satisfy the Fibered Isomorphism Conjecture if $G_1$ and $G_2$ are
  virtually cyclic. Since every virtually cyclic group admits an
  epimorphism to $\IZ$ or $D_{\infty}$ with finite kernel and $\IZ
  \subseteq D_{\infty}$, the product $G_1 \times G_2$ admits a group
  homomorphism to $D_{\infty} \times D_{\infty}$ with finite kernel.
  Now apply Lemma~\ref{lem:Fibered_Conjecture_passes_to_subgroups} and
  Lemma~\ref{lem:extensions_with_finite_kernel}.
\end{proof}

Assertion~\ref{lem:products:closed_under_products} appearing in
Lemma~\ref{lem:products} is interesting in the case $\calc = \Fin$.

\begin{lemma} \label{lem:elementary_amenable}
  Suppose that $\calh^?_*$ is continuous.  Suppose that any virtually
  finitely generated abelian group satisfies the Fibered Isomorphism
  Conjecture for $\Fin$.

  Then every elementary amenable group satisfies the Fibered
  Isomorphism Conjecture for $\Fin$.
\end{lemma}
\begin{proof}
  Using the same transfinite induction strategy and the same notation
  as in the proof of \cite[Corollary~3.9]{Farrell-Linnell(2003b)}, one
  reduces the claim to the following assertion. The group $G$
  satisfies the Fibered Isomorphism Conjecture for $\Fin$ provided
  that there exists an extension $1 \to H \to G \xrightarrow{p} A \to
  1$ such that $A$ is virtually finitely generated abelian, $H$ belongs to
  $L\calx_{\alpha-1}$ and the Fibered Isomorphism Conjecture for
  $\Fin$ holds for every group in the class of groups
  $\calx_{\alpha-1}$. Here $L\calx_{\alpha-1}$ is the class of groups
  for which every finitely generated subgroup occurs as a subgroup of
  some group in $\calx_{\alpha-1}$.
  Theorem~\ref{the:isomorphism_conjecture_is_stable_under_colim}~\ref{the:isomorphism_conjecture_is_stable_under_colim:injective}
  and Lemma~\ref{lem:Fibered_Conjecture_passes_to_subgroups} 
  imply
  that the Fibered Isomorphism Conjecture holds for 
  every group in $L\calx_{\alpha-1}$.
  Because of
  Lemma~\ref{lem:extensions_fibered_conclusion} it remains to prove
  for any finite subgroup $K \subseteq A$ that the Fibered Isomorphism
  Conjecture for $\Fin$ holds for $p^{-1}(K)$.
  We conclude from the short exact sequence $1 \to H \to p^{-1}(K) \to K \to 1$
  and \cite[Lemma~2.1~(iii)]{Farrell-Linnell(2003b)} that $p^{-1}(K)$ is a
  member of $L\calx_{\alpha-1}$ and satisfies therefore the Fibered Isomorphism Conjecture
  for $\Fin$.
\end{proof}

\begin{lemma} \label{lem:virtually_nilpotent}
  Suppose that $\calh^?_*$ is continuous.  Suppose that any virtually
  finitely generated abelian group satisfies the Fibered Isomorphism
  Conjecture for $\VCyc$.

  Then every virtually nilpotent group satisfies the Fibered
  Isomorphism Conjecture for $\VCyc$.
\end{lemma}

\begin{proof}
  Any finitely generated subgroup of a virtually abelian group is
  virtually finitely generated abelian.
  The assumptions and  Theorem~\ref{the:isomorphism_conjecture_is_stable_under_colim}~\ref{the:isomorphism_conjecture_is_stable_under_colim:injective}
  imply that any virtually abelian group satisfies
  the Fibered Isomorphism Conjecture for $\VCyc$.

  For a group $H$ we denote by $\cent(H)$ its center.
  Recall that a group $N$ is called \emph{nilpotent}, if we can find a
  sequence of epimorphisms $N = N_0 \xrightarrow{p_0} N_1
  \xrightarrow{p_1} \cdots \xrightarrow{p_r} N_r$ such that $\ker(p_i)
  = \cent(N_{i-1})$ for $i = 1,2, \ldots, r$ and $N_r = \{1\}$. The
  class of $N$ is the smallest non-negative integer $r$ for which such
  a sequence of epimorphisms exists. Let $G$ be virtually nilpotent.
  Hence we can find a normal subgroup $N \subseteq G$ such that $N$ is
  nilpotent and $G/N$ is finite.  We show by induction over the class
  of $N$ that $G$ satisfies the Fibered Isomorphism Conjecture for
  $(G,\VCyc)$. The induction beginning, where the class is $\le 1$ and
  hence $G$ is virtually abelian, has already been taken care of.

  We can arrange that $r$ is the class of $N$.
  Since $N$ is normal in
  $G$ and $\cent(N)$ is a characteristic subgroup of $N$, $\cent(N)$
  is a normal subgroup of $G$. We obtain the exact sequence $1 \to
  N/\cent(N) \to G/\cent(N) \to G/N \to 1$. Since the class of
  $N/\cent(N)$ is smaller than the class of $N$, the Fibered
  Isomorphism Conjecture holds for $(G/\cent(N),\VCyc)$ by the
  induction hypothesis. Because of
  Theorem~\ref{lem:extensions_fibered_conclusion} it remains to show
  for any virtually cyclic subgroup $V \subseteq G/\cent(N)$ that
  $(q^{-1}(V),\VCyc)$ satisfies the Fibered Isomorphism Conjecture,
  where $q \colon G \to G/\cent(N)$ is the canonical projection.

  Let $\phi \colon G \to \aut(N)$ be the group homomorphism sending $g
  \in G$ to the automorphism of $N$ given by conjugation with $g$.
  Since $\cent(N)$ is a characteristic subgroup of $N$, it induces a
  homomorphism $\phi' \colon G \to \aut(\cent(N))$. Since the
  conjugation action of $N$ on itself is the identity on $\cent(N)$,
  the homomorphism $\phi'$ factorizes through the finite group $G/N$
  and hence has finite image. Hence we can find $g \in G$ such that
  $q(g)$ generates an infinite cyclic subgroup $C$ in $V$ whose index
  in $V$ is finite and $\phi'(g) = \id_{\cent(N)}$.  Hence $q^{-1}(C)$
  is isomorphic to $\cent(N) \times C$ and has finite index in
  $q^{-1}(V)$. Therefore $q^{-1}(V)$ is virtually abelian and
  $(q^{-1}(V),\VCyc)$ satisfies the Fibered Isomorphism Conjecture.
\end{proof}


\section{The Farrell-Jones Conjecture for algebraic $K$-Theory}
\label{sec:The_Farrell-Jones_Conjecture_for_algebraic_K-Theory}

Recall that the \emph{(Fibered) Farrell-Jones Conjecture for algebraic
$K$-theory with coefficients in $R$ for the group $G$} is the
(Fibered) Isomorphism Conjecture~\ref{def:(Fibered)_Isomorphism_Conjectures_for_calh?_ast}
in the special case, where the family
$\calf$ consists of all virtually cyclic subgroups of $G$ and
$\calh^?_*$ is the equivariant homology theory $\calh^?_*(-;\bfK_R)$
associated to the $\Groupoids$-spectrum given by algebraic $K$-theory
and $R$ as coefficient ring (see \cite[Section~6]{Lueck-Reich(2005)}).
So the Farrell-Jones Conjecture for algebraic
$K$-theory with coefficients in $R$ for the group $G$ predicts that
the map
$$H_n^G(\EGF{G}{\VCyc},\bfK_R) \to K_n(RG)$$
is bijective for all $n \in \IZ$. The original source for
(Fibered) Farrell-Jones Conjecture is~\cite[1.6 on page~257 and~1.7 on page~262]{Farrell-Jones(1993a)}.

Recall that the \emph{(Fibered) Farrell-Jones Conjecture for homotopy
$K$-theory with coefficients in $R$ for the group $G$} is the
(Fibered) Isomorphism Conjecture~\ref{def:(Fibered)_Isomorphism_Conjectures_for_calh?_ast}
in the special case, where the family
$\calf$ consists of all finite subgroups of $G$ and $\calh^?_*$ is the
equivariant homology theory $\calh^?_*(-;\bfKH_R)$ associated to the
$\Groupoids$-spectrum given by homotopy $K$-theory and $R$ as
coefficient ring (see \cite[Section~7]{Bartels-Lueck(2006)}).  So the
Farrell-Jones Conjecture for homotopy $K$-theory with coefficients in
$R$ for the group $G$ predicts that the map
$$H_n^G(\EGF{G}{\Fin},\bfKH_R) \to \KH_n(RG)$$
is bijective for all $n
\in \IZ$.

The following theorem follows from the main result of
Bartels-L\"uck-Reich~\cite{Bartels-Lueck-Reich(2007hyper)}
together with~\cite[Corollary~4.3]{Bartels-Reich(2005)}.

\begin{theorem}
\label{the:FiberedFJC_for_word-hyperbolic_groups} The Fibered
Farrell-Jones Conjecture for algebraic $K$-theory is true for every
word-hyperbolic group and every coefficient ring.
  \end{theorem}

  Next we extend a result of Quinn~\cite[Theorem~1.2.2]{Quinn(2005)}
  for virtually abelian groups to virtually nilpotent groups.

 \begin{theorem}
  \label{the:FiberedFJC_for_virtually_nilpotent_groups}
  The Fibered Farrell-Jones Conjecture for algebraic $K$-theory is
  true for every virtually nilpotent group and every coefficient ring.
  \end{theorem}
  \begin{proof}
    The Fibered Farrell-Jones Conjecture for algebraic $K$-theory is
    true for every virtually abelian group and every coefficient ring
    by a result of Quinn~\cite[Theorem~1.2.2]{Quinn(2005)}.  (Quinn
    deals only with commutative coefficient rings but this assumption
    is not needed in his argument.) Now apply
    Lemma~\ref{lem:virtually_nilpotent}.
  \end{proof}

  For more information about groups satisfying the Farrell-Jones
  Conjecture for algebraic $K$-theory with coefficients in $\IZ$
  and the Farrell-Jones Conjecture for homotopy $K$-theory
  $\KH_*(RG)$ we refer to~\cite[Sections~5.2 and~5.3]{Lueck-Reich(2005)}
  and~\cite[Theorem~0.5]{Bartels-Lueck(2006)}.

\begin{lemma} \label{lem:H(-,K_R)_strongly_continuous)}
  The equivariant homology theories $\calh^?_*(-;\bfK_R)$ and
  $\calh^?_*(-;\bfKH_R)$ are strongly continuous.
\end{lemma}
\begin{proof}
  We have to show for every directed systems of groups $\{G_i \mid i \in
  I\}$ with $G = \colim_{i \in I} G_i$ that the canonical maps
  \begin{eqnarray*}
  \colim_{i \in I} K_n(RG_i) & \to & K_n(RG);
  \\
  \colim_{i \in I} \KH_n(RG_i) & \to & \KH_n(RG),
  \end{eqnarray*}
  are bijective for all $n \in \IZ$.  Obviously $RG$ is the colimit of
  rings $\colim_{i \in I} RG_i$. Now the claim follows for $K_n(RG)$ for
  $n \ge 0$ from \cite[(12) on page~20]{Quillen(1973)}.  Using the
  Bass-Heller-Swan decomposition one gets the results for $K_n(RG)$
  and also for the Nil-groups $N^pK_n(RG)$ defined by
  Bass~\cite[XII]{Bass(1968)} for all $n \in \IZ$ and $p \ge 1$.  The
  claim for $\KH_n(RG)$ follows from the spectral sequence
  \cite[Theorem~1.3]{Weibel(1989)}.
\end{proof}

We conclude from Theorem~\ref{the:isomorphism_conjecture_is_stable_under_colim}
and Lemma~\ref{lem:H(-,K_R)_strongly_continuous)}
that the (Fibered) Farrell-Jones Conjecture is inherited under directed colimits.


\subsection{Extensions}
\label{subsec:Extensions}

\begin{lemma} \label{lem:extensions_with_kernel_virtually_cyclic}
  Let $1 \to K \to G \to Q \to 1$ be an extension of groups.  Suppose
  that $K$ is virtually cyclic and $Q$ satisfies the Fibered Farrell-Jones
  Conjecture with coefficients in $R$.

  Then $G$ satisfies the Fibered Farrell-Jones Conjecture with
  coefficients in $R$.
\end{lemma}
\begin{proof}
  Because of Lemma~\ref{lem:extensions_fibered_conclusion} it suffices
  to prove that $(G;\VCyc(G))$ satisfies the Fibered Farrell Jones
  Conjecture in the case that $Q$ is virtually cyclic. Choose an
  infinite cyclic subgroup $C$ of $Q$. Let $\phi \colon K \to K$ be
  the automorphism given by conjugation with an element in $G$ which
  is mapped to a generator of $C$ under the epimorphism $G \to Q$.
  Then $p^{-1}(C)$ is a subgroup of $G$ which has finite index and is
  isomorphic to the semi-direct product $K \rtimes_{\phi} C$.  Since
  $K$ is virtually cyclic, its automorphism group has finite order.
  Hence by replacing $C$ by a subgroup of the order of this
  automorphism group as index, we can arrange that $p^{-1}(C)$ is a
  subgroup of finite index in $G$ and $p^{-1}(C) \cong K \times C$.
  Since $K$ is virtually cyclic, we conclude that $G$ contains a
  subgroup of finite index which is isomorphic to $\IZ^2$.  In
  particular $G$ is virtually abelian. Hence $(G,\VCyc(G))$ satisfies
  the Fibered Farrell-Jones Conjecture by
  Theorem~\ref{the:FiberedFJC_for_virtually_nilpotent_groups}.
\end{proof}

\begin{lemma} \label{lem:FFJC_and_products}
  Let $G_1$ and $G_2$ be groups.  Then $G_1 \times G_2$ satisfies the
  Fibered Farrell-Jones Conjecture with coefficients in $R$ if and
  only if both $G_1$ and $G_2$ satisfy the Fibered Farrell-Jones
  Conjecture with coefficients in $R$.
\end{lemma}
\begin{proof}
  Because of Lemma~\ref{lem:products}~\ref{lem:products:VCcy} it
  suffices to show that $D_{\infty} \times D_{\infty}$ satisfies the
  Fibered Farrell-Jones Conjecture.  This follows from
  Theorem~\ref{the:FiberedFJC_for_virtually_nilpotent_groups}.
\end{proof}

Next we give the proof of Theorem~\ref{the:class_FJC}.

\begin{proof}\ref{the:class_FJC:word_hyper_virt_nilpotent} This follows from
  Theorem~\ref{the:FiberedFJC_for_word-hyperbolic_groups} and
  Theorem~\ref{the:FiberedFJC_for_virtually_nilpotent_groups}.
  \\[1mm]\ref{the:class_FJC:products} This follows from
  Theorem~\ref{lem:FFJC_and_products}.
  \\[1mm]\ref{the:class_FJC:colim} This follows from
  Theorem~\ref{the:isomorphism_conjecture_is_stable_under_colim}~\ref{the:isomorphism_conjecture_is_stable_under_colim:general}
  and
  Lemma~\ref{lem:H(-,K_R)_strongly_continuous)}.
  \\[1mm]\ref{the:class_FJC:subgroups} This follows from
  Lemma~\ref{lem:Fibered_Conjecture_passes_to_subgroups}.
\end{proof}

A ring $R$ is called \emph{regular coherent} if every finitely
presented $R$-module possesses a finite-dimensional resolution by
finitely generated projective $R$-modules.  A ring $R$ is
\emph{regular} if and only if it is regular coherent and Noetherian.
A group $G$ is called \emph{regular} or \emph{regular coherent}
respectively if for any regular ring $R$ the group ring $RG$ is
regular respectively regular coherent.  Poly-$\IZ$-groups and free
groups are regular coherent (see~\cite[page~247]{Waldhausen(1978a)}.
For more information about these notions we refer to \cite[Theorem
19.1]{Waldhausen(1978a)}.

\begin{theorem} \label{the:torsionfree_extension_CL-hyperbolic}
  Suppose that $Q$ is torsionfree and that the Fibered Farrell-Jones
  Conjecture holds for $Q$.  Suppose that $K$ is a regular coherent
  group. Suppose that $R$ is regular. Let $1 \to K \xrightarrow{i}
  G \xrightarrow{p} Q \to 1$ be an extension of groups.

  Then the assembly map
  $$H_n(BG;\bfK_R) \xrightarrow{\cong} K_n(RG)$$
  is an isomorphism for
  $n \in \IZ$.
\end{theorem}

\begin{proof}
  Let $\Gamma$ be a torsionfree group. The relative assembly map
  $H_n(B\Gamma;\bfK_R) \xrightarrow{\cong}
  H_n^{\Gamma}(\EGF{\Gamma}{\VCyc};K_{R})$ is bijective
  (see~\cite[Proposition~2.6 on page~686]{Lueck-Reich(2005)}) since
  $R$ is regular and $\Gamma$ is torsionfree.  Hence the Farrell-Jones
  Conjecture for $\Gamma$ and $R$ boils down to the claim that the
  assembly map $H_n(B\Gamma;\bfK_R) \xrightarrow{\cong} K_n(R\Gamma)$
  is bijective for all $n \in \IZ$. This implies that
  Theorem~\ref{the:torsionfree_extension_CL-hyperbolic} follows
  directly from Lemma~\ref{lem:extensions_unfibered_conclusion} if we
  can show that the second assumption appearing in
  Lemma~\ref{lem:extensions_unfibered_conclusion} is satisfied.

  Let $\calcl$ be the class of groups introduced in~\cite[Definion~19.2 on page
  248]{Waldhausen(1978a)} or~\cite[Definition~0.10]{Bartels-Lueck(2006)}.  Let $V
  \subseteq Q$ be virtually cyclic. Since $K$ is regular coherent and $Q$ is torsionfree
  and hence $V$ is isomorphic to $\IZ$, $p^{-1}(V)$ belongs to the class $\calcl$. Since
  $\calcl$ by \cite[Proposition~19.3 on page~249]{Waldhausen(1978a)} is closed under taking subgroups,
  every element in $p^*\VCyc(Q)$ belongs to $\calcl$. One of
  the main results in Waldhausen's article \cite{Waldhausen(1978a)} is that for a regular
  ring $R$ the $K$-theoretic assembly map $H_n(BG';\bfK_R) \to K_n(RG')$ is an isomorphism
  for $G' \in \calcl$. Hence the second assumption appearing in
  Lemma~\ref{lem:extensions_unfibered_conclusion} is satisfied.
\end{proof}


\subsection{Passage from $\Fin$ to $\VCyc$}
\label{subsec:Passage_from_Fin_to_VCyc}

\begin{lemma} \label{lem:relative_ass_Fin_to_Vcycv_iso_char_zero}
  Let $G$ be a group. Let $R$ be a regular ring such that the order of
  any finite subgroup of $G$ is invertible in $R$.

  Then the relative assembly map
  $$
  H_n^G ( \EGF{G}{\Fin} ; \bfK_R ) \to H_n^G ( \EGF{G}{\VCyc} ; \bfK_R)
  $$
  (induced by the up to $G$-homotopy unique $G$-map $\EGF{G}{\Fin} \to \EGF{G}{\VCyc}$)
  is an isomorphism for all $n \in \IZ$.
\end{lemma}

\begin{proof}
This is proved for instance in~\cite[Proposition~2.6 on page~686]{Lueck-Reich(2005)}.
\end{proof}

\begin{lemma} \label{lem:relative_ass_Fin_to_Vcycv_R_is_Z}
  Let $G$ be a group. Then the relative assembly map induces for all
  $n \in \IZ$ isomorphisms
  $$
  H_n^G ( \EGF{G}{\Fin}  ; \bfK_{\IZ}) \otimes_{\IZ} \IQ \to 
  H_n^G ( \EGF{G}{\VCyc} ; \bfK_{\IZ}) \otimes_{\IZ} \IQ.$$
\end{lemma}
\begin{proof} This is proved  in~\cite[Theorem~5.6]{Grunewald(2006)}.
\end{proof}

Notice that the two
Lemmas~\ref{lem:relative_ass_Fin_to_Vcycv_iso_char_zero}
and~\ref{lem:relative_ass_Fin_to_Vcycv_R_is_Z} above deal
\emph{not} with the Fibered version.  A discussion of the Fibered version of
Lemma~\ref{lem:relative_ass_Fin_to_Vcycv_iso_char_zero} can be
found in Subsection~\ref{subsec:Nil-groups_and_regular_rings_R_with_Z_subseteq_R}.
The fibered version of
Lemma~\ref{lem:relative_ass_Fin_to_Vcycv_R_is_Z} is definitely false.
For homotopy $K$-theory one can prove in the fibered situation that
the passage from $\Fin$ to $\VCyc$ does not matter.
For $\Lambda = \IZ$ the next lemma
is proven in~\cite[Remark~7.4]{Bartels-Lueck(2006)}.
The following more general statement follows by the same argument.

\begin{lemma} \label{lem:FJC_fibered_Fin_Vcyc_homotopy_K-theory}
  Let $R$ be a ring. A group $G$ satisfies the (Fibered) Isomorphism
  Conjecture for homotopy $K$-theory with coefficients in $R$ for the
  family $\Fin$ after applying $- \otimes_{\IZ} \Lambda$ for $\IZ
  \subseteq \Lambda \subseteq \IQ$ to the assembly map
  if and only if $G$ satisfies the (Fibered) Isomorphism
  Conjecture for homotopy $K$-theory with coefficients in $R$ for the
  family $\VCyc$ after applying $- \otimes_{\IZ} \Lambda$ for $\IZ
  \subseteq \Lambda \subseteq \IQ$ to the assembly map.
\end{lemma}


\subsection{Homotopy $K$-theory and rings with finite characteristic}
\label{subsec:Homotopy_K-theory_and_rings_with_finite_characteristic}

\begin{lemma} \label{lem:Vanishing_of_NK_p(R)[1/N]}
  Let $R$ be a ring of finite characteristic $N$.
  Let $N\!K_n(R)$ be the Nil-group of Bass.
  Then we get for $n \in \IZ$
  $$N\!K_n(R)[1/N]~=~0.$$
\end{lemma}
\begin{proof}
  The proof can be found in~\cite[Corollary~3.2]{Weibel(1981)}. We
  give a brief outline for the reader's convenience.

  Put $\Lambda = \IZ/N$. Then $R$ is a $\Lambda$-algebra. Let
  $W(\Lambda)$ be the ring of big Witt vectors over $\Lambda$. The
  underlying additive group is the multiplicative group $1 + t \Lambda[[t]]$
  of formal power series with leading term $1$. We do not need the
  explicit multiplicative structure but need to know that the identity element is
  $1 -t$. Let $\End (\Lambda)$ be the Grothendieck group of
  endomorphisms $f \colon P \to P$ of finitely generated projective
  $\Lambda$-modules.  We get an injective homomorphism $K_0(\Lambda)
  \to \End (\Lambda)$ by sending $[P]$ to $[0 \colon P \to P]$. Its
  cokernel is denoted by $\End_0(\Lambda)$.  The tensor product
  induces the structure of a commutative ring on $\End(\Lambda)$ for
  which $K_0(\Lambda )$ becomes an ideal. Hence $\End_0(\Lambda)$ is a
  commutative ring.  There is a $\End_0(\Lambda)$-module structure on
  $N\!K_n(R)$ for all $n \in \IZ$.  Almkvist~\cite{Almkvist(1973)} shows
  that the characteristic polynomial defines an injective ring
  homomorphism
  $$\chi \colon \End_0(\Lambda) \to W(\Lambda), \quad [f\colon P \to P]
  \mapsto \det(\id_P - t \cdot f).$$
  For a positive integer $N$ let
  $I_N$ be the ideal in $\End_0(\Lambda)$
  $$I_N = \{[f \colon P \to P] \mid \chi([f]) \equiv 1 \mod t^N\}.$$
  Stienstra~\cite{Stienstra(1982)} (see also Theorem~1.3 and the
  following paragraph in \cite{Weibel(1981)}) has
  proven that for every element $x \in N\!K_n(R)$ there exists $N(x)$
  such that the ideal $I_{N(x)}$ annihilates $x$.  Now choose a
  sufficiently large positive integer $k$ such that $\binom{N^k}{j}$
  is a multiple of $N$ for $1 \le j < N(x)$. Then we conclude
  $$\chi(N^k \cdot [\id_{\Lambda}]) = (1-t)^{N^k} = \sum_{j = 0}^{N^k}
  \binom{N^k}{j} (-t)^j = 1 + \sum_{j = N(x) }^{N^k} \binom{N^k}{j} (-t)^j
  \equiv 1 \mod t^{N(x)}.$$
  This shows that $N^k \cdot [\id_{\Lambda}] \in
  I_{N(x)}$. We compute
  $$N^k \cdot x = N^k \cdot \left([\id_{\Lambda}] \cdot x\right) =
  \left(N^k \cdot [\id_{\Lambda}]\right) \cdot x = 0.$$
  This implies
  $N\!K_n(R)[1/N]=0$.
\end{proof}

\begin{lemma} \label{lem:K-to-KH}
\begin{enumerate}
\item \label{lem:K-to-KH:finite}
  Let $R$ be a ring of finite characteristic $N$.
  Then the canonical map from algebraic $K$-theory to homotopy $K$-theory
  induces an isomorphism
  $$K_n(R)[1/N] \xrightarrow{\cong} \KH_n(R)[1/N]$$
  for all $n \in \IZ$;

\item \label{lem:K-to-KH:IZ}
  Let $H$ be a finite group.
  Then the canonical map from algebraic $K$-theory to homotopy $K$-theory
  induces an isomorphism
  $$K_n( \IZ[H] ) \otimes_\IZ \IQ  \xrightarrow{\cong} \KH_n( \IZ[H] ) \otimes_\IZ \IQ$$
  for all $n \in \IZ$.
\end{enumerate}
\end{lemma}

\begin{proof}\ref{lem:K-to-KH:finite}
  We conclude from Lemma~\ref{lem:Vanishing_of_NK_p(R)[1/N]} that
  $N\!K_n(R)[1/N] = 0$. This implies that $N^pK_n(R)[1/N] = 0$ for $n
  \in \IZ$ and $p \ge 1$. Now apply the spectral sequence from
  \cite[Theorem~1.3]{Weibel(1989)}.
  \\[1mm]\ref{lem:K-to-KH:IZ}
  By \cite[Corollary]{Weibel(1981)} and \cite[Remark~8.3]{Bartels-Lueck(2006)}
  $N\!K_n(\IZ[H]) \otimes_\IZ \IQ = 0$.
  Now proceed as in~\ref{lem:K-to-KH:finite}.
\end{proof}

\begin{lemma} \label{lem:assembly-for-K-versus-KH}
Let $(X,A)$ be a pair of $G$-$CW$-complexes.
\begin{enumerate}
\item \label{lem:assembly-for-K-versus-KH:finite}
Let $R$ be a ring of finite characteristic $N$.
Then the natural map
$$H_n^G(X,A;\bfK_R)[1/N] \to H_n^G(X,A;\bfKH_R)[1/N]$$
is bijective for every $n \in \IZ$;

\item \label{lem:assembly-for-K-versus-KH:IZ}
Assume that $X$ is a $G$-$\Fin$-$CW$-complex, i.e.,
the isotropy groups of $X$ are finite.
Then the natural map
$$H_n^G(X,A;\bfK_\IZ) \otimes_\IZ \IQ \to H_n^G(X,A;\bfKH_\IZ) \otimes_\IZ \IQ$$
is bijective for every $n \in \IZ$.
\end{enumerate}
\end{lemma}

\begin{proof}
  This follows from Lemma~\ref{lem:K-to-KH} and a spectral
  sequence argument based on the equivariant Atiyah-Hirzebruch spectral
  sequence~(see for instance~\cite[Theorem~4.7]{Davis-Lueck(1998)}).
\end{proof}

\begin{lemma}
\label{lem:FJC_for_K_versus_KH}
Let $R$ be a ring and let $G$ be a group. Let $N \ge 2$ be an integer.
\begin{enumerate}

\item \label{lem:FJC_for_K_versus_KH:K_implies_KH} If $G$ satisfies
      the (Fibered) Farrell-Jones Conjecture for algebraic $K$-theory
      with coefficients in $R[x_1,x_2, \ldots , x_k]$ for all $k \ge 0$,
      then $G$ satisfies the (Fibered) Farrell-Jones Conjecture for
      homotopy $K$-theory with coefficients in $R$;

\item \label{lem:FJC_for_K_versus_KH:KH_implies_K} Suppose that $N
      \cdot 1_R = 0$ and that $G$ satisfies the (Fibered)
      Farrell-Jones Conjecture for homotopy $K$-theory with
      coefficients in $R$ after applying $-\otimes_{\IZ} \IZ[1/N]$.
      Then $G$ satisfies the (Fibered) Farrell-Jones Conjecture for
      algebraic $K$-theory with coefficients in $R$ after applying
      $-\otimes_{\IZ} \IZ[1/N]$ for both the family $\Fin$ and
      $\VCyc$.
\end{enumerate}
\end{lemma}

\begin{proof}\ref{lem:FJC_for_K_versus_KH:K_implies_KH} This is proven in
  \cite[Theorem~8.4]{Bartels-Lueck(2006)}.
  \\[2mm]\ref{lem:FJC_for_K_versus_KH:KH_implies_K}
  Consider for any family $\calf$ of subgroups of $G$ the following commutative diagram
  $$
  \xycomsquare{H_n^G(\EGF{G}{\calf};\bfK_R)[1/N]}{}{H_n(\pt;\bfK_R)[1/N]
    = K_n(RG))[1/N]} {\cong}{\cong}
  {H_n^G(\EGF{G}{\calf};\bfKH_R)[1/N]}{}{H_n(\pt;\bfKH_R)[1/N] =
    \KH_n(RG)[1/N]}
  $$
  where the horizontal maps are the assembly maps induced by the
  projection $\EGF{G}{\calf} \to \pt$ and the vertical maps are
  induced by the passage from algebraic $K$-theory to homotopy $K$-theory.
  Lemma~\ref{lem:assembly-for-K-versus-KH}~\ref{lem:assembly-for-K-versus-KH:finite} implies
  that the vertical maps are bijective.
  Now apply Lemma~\ref{lem:FJC_fibered_Fin_Vcyc_homotopy_K-theory}.
\end{proof}

Next we can give the proof of Theorem~\ref{the:class_FJ_KH_andFJ_N}.

\begin{proof}\ref{the:class_FJ_KH_andFJ_N:elem_amen_word_hyper} Word-hyperbolic
  groups and virtually abelian groups satisfy the Farrell-Jones
  Conjecture for algebraic $K$-theory with coefficients in any ring
  $R$ by Theorem~\ref{the:FiberedFJC_for_word-hyperbolic_groups} and
  Theorem~\ref{the:FiberedFJC_for_virtually_nilpotent_groups}.  We
  conclude from
  Lemma~\ref{lem:FJC_for_K_versus_KH}~\ref{lem:FJC_for_K_versus_KH:K_implies_KH}
  that word-hyperbolic groups and virtually abelian groups belong to
  $\calfj_{\KH}(R)$ for all rings $R$.  We conclude from
  Lemma~\ref{lem:elementary_amenable} that all elementary amenable groups
  belong to $\calfj_{\KH}(R)$ for all rings $R$.
  Lemma~\ref{lem:FJC_for_K_versus_KH}~\ref{lem:FJC_for_K_versus_KH:KH_implies_K}
  implies that all word-hyperbolic groups and all elementary amenable
  groups belong to $\calfj_{FC}(R)$.

  In particular every virtually cyclic group belongs to $\calfj_{FC}(R)$. We conclude from
  Theorem~\ref{the:transitivity} that for a ring $R$ of characteristic $N$ a group
  satisfies the Fibered Farrell-Jones Conjecture for algebraic $K$-theory for $G$ with
  coefficients in $R$ after applying $- \otimes_{\IZ} \IZ[1/N]$ to the assembly map for
  the family $\Fin$ if and only it does for the family $\VCyc$.  Hence in the sequel we
  only have to consider the family $\Fin$ when dealing with $\calfj_{FC}(R)$.
  \\[1mm]\ref{the:class_FJ_KH_andFJ_N:products} This follows from
  Lemma~\ref{lem:products}~\ref{lem:products:closed_under_products}.
  \\[1mm]\ref{the:class_FJ_KH_andFJ_N:colim} This follows from
  Theorem~\ref{the:isomorphism_conjecture_is_stable_under_colim}~\ref{the:isomorphism_conjecture_is_stable_under_colim:general} 
  and Lemma~\ref{lem:H(-,K_R)_strongly_continuous)}. 
  \\[1mm]\ref{the:class_FJ_KH_andFJ_N:subgroups} This follows from
  Lemma~\ref{lem:Fibered_Conjecture_passes_to_subgroups}.
  \\[1mm]\ref{the:class_FJ_KH_andFJ_N:extensions} This follows from
  Lemma~\ref{lem:extensions_fibered_conclusion} since a group which is commensurable to a
  word-hyperbolic group is again word-hyperbolic and the same is true for elementary
  amenable groups.
  \\[1mm]\ref{the:class_FJ_KH_andFJ_N:actions_on_trees} This is proven for $\calfj_{\KH}$ in
  Bartels-L\"uck~\cite[Theorem~0.5]{Bartels-Lueck(2004ind)}. The same proof applies to
  $\calfj_{FC}$ if we can show that $\calh^?_*(-;\bfK_R)[1/N]$ has the tree property
  (see~\cite[Definition~4.1 and Theorem~4.2]{Bartels-Lueck(2004ind)}). This follows from
  Lemma~\ref{lem:assembly-for-K-versus-KH}~\ref{lem:assembly-for-K-versus-KH:finite} since the equivariant
  homology theory $\calh^?_*(-;\bfKH_R)$ has the tree property by
  \cite[Theorem~11.1]{Bartels-Lueck(2004ind)}.

  The claim $\calfj_{\KH}(R) \subseteq \calfj_{FC}(R)$ follows from
  Lemma~\ref{lem:FJC_for_K_versus_KH}~\ref{lem:FJC_for_K_versus_KH:KH_implies_K}.

\end{proof}


\subsection{Coefficient rings with operation}
\label{subsec:Coefficent_rings_with_operation}

In the setup developed so far we have not dealt with the more general
version developed in Bartels-Reich~\cite{Bartels-Reich(2005)} where
one fixes a group and considers an additive category with $G$-action. This setup
can deal with crossed products $R \rtimes G$ and not only with group
rings $RG$.  However, a slight modification of the
proofs above allows to carry over the result above to this setting.
This is explained in Bartels-Echterhoff-L\"uck~\cite{Bartels-Echterhoff-Lueck(2007colim)}.


\section{The projective class group and induction from finite subgroups}
\label{sec:The_projective_class_group_and_induction_from_finite_subgroups}

Let $\OrGF{G}{\Fin}$ be the category whose objects are homogeneous
spaces $G/H$ with finite $H$ and whose morphisms are $G$-maps. We
obtain a functor from $\OrGF{G}{\Fin}$ to abelian groups by sending
$G/H$ to $K_0(RH)$. It sends a morphism $G/H \to G/K, gH \mapsto
gg_0K$ to the map $K_0(RH) \to K_0(RK)$ coming from the group
homomorphism $H \to K, h \mapsto g_0hg_0^{-1}$.  This is well-defined
since inner automorphisms of $H$ induce the identity on $K_0(RH)$.
The various inclusions of finite subgroups of $G$ yield a homomorphism
\begin{eqnarray}
I(G;R) \colon \colim_{\OrGF{G}{\Fin}} K_0(RH) & \to & K_0(RG).
\label{I(G,R):colim_FIN_to_K_0(RG)}
\end{eqnarray}
Notice for the sequel that the canonical map of $\Lambda$-modules
$$\left(\colim_{\OrGF{G}{\Fin}} K_0(RH)\right) \otimes_{\IZ} \Lambda
\xrightarrow{ \cong} \colim_{\OrGF{G}{\Fin}} \left(K_0(RH)
  \otimes_{\IZ} \Lambda\right)$$
is bijective for every ring $\Lambda$.

The next lemma is proven in \cite[Lemma~2.9]{Farrell-Linnell(2003b)}
for fields and carries over directly to skew-fields.

\begin{lemma} \label{lem:Negative_K-groups_of_FH-Vanish}
  Let $D$ be a skew-field (of arbitrary characteristic) and $H$ be a finite
  group. Then $K_n(DH) = 0$ for $n \le -1$.
\end{lemma}

Now we can give the proof of
Theorem~\ref{the:FJC_implies_Moodys_induction}
\begin{proof}\ref{the:FJC_implies_Moodys_induction:reg_ring_Q_subset_R} This is
  proved in~\cite[page~691]{Lueck-Reich(2005)}.
  \\[1mm]\ref{the:FJC_implies_Moodys_induction:Field_char_p} This follows
  analogously to the proof of
  assertion~\ref{the:FJC_implies_Moodys_induction:reg_ring_Q_subset_R}
  using Lemma~\ref{lem:Vanishing_of_NK_p(R)[1/N]} and
  Lemma~\ref{lem:Negative_K-groups_of_FH-Vanish}.
\end{proof}

If $G$ satisfies the Farrell-Jones Conjecture for algebraic $K$-theory
with coefficients in $\IZ$, then the following maps are injective
(see~\cite[page~692]{Lueck-Reich(2005)}).
\begin{eqnarray*}
\colim_{\OrGF{G}{\Fin}} K_n( \IZ H) \otimes_{\IZ} \IQ & \to & K_n( \IZ G)
\otimes_{\IZ} \IQ \quad \text{ for }n \in \IZ;
\\
\colim_{\OrGF{G}{\Fin}} \Wh(H) \otimes_{\IZ} \IQ & \to & \Wh(G) \otimes_{\IZ} \IQ.
\end{eqnarray*}
The injectivity for the map involving the Whitehead group is proven
for groups $G$ satisfying a  mild homological finiteness conditions in
\cite {Lueck-Reich-Rognes-Varisco(2007)}.  In general these maps are
not surjective.  In particular
$$\colim_{\OrGF{G}{\Fin}} K_0(\IZ H) \otimes_{\IZ} \IQ \to K_n(\IZ G)
\otimes_{\IZ} \IQ$$
is in general not surjective.


\section{Bass Conjectures}
\label{sec:Bass_Conjectures}

In this section we explain the relationship between the Farrell-Jones
Conjecture for algebraic $K$-theory and the Bass Conjecture.

Let $G$ be a group.  Let $\con(G)$ be the set of conjugacy classes
$(g)$ of elements $g \in G$. Denote by $\con(G)_f$ the subset of
$\con(G)$ consisting of those conjugacy classes $(g)$ for which each
representative $g$ has finite order. Let $R$ be a commutative ring.
Let $\class_R(G)$ and $\class_R(G)_f$ be the free $R$-module with the
set $\con(G)$ and $\con(G)_f$ as basis.  This is the same as the
$R$-module of $R$-valued functions on $\con(G)$ and $\con(G)_f$ with
finite support.
Define the \emph{universal $R$-trace}%
\begin{eqnarray}  \tr_{RG}^u \colon RG \to \class_R(G), \quad
\sum_{g \in G} r_g \cdot g~\mapsto~\sum_{g \in G} r_g \cdot (g).
\label{universal_CG-trace}
\end{eqnarray}
It extends to a function $\tr_{RG}^u \colon M_n(RG) \to \class_R(G)$
on $(n,n)$-matrices over $RG$ by taking the sum of the traces of the
diagonal entries.  Let $P$ be a finitely generated projective
$RG$-module. Choose a matrix $A \in M_n(RG)$ such that $A^2 = A$ and
the image of the $RG$-map $r_A \colon RG^n \to RG^n$ given by right
multiplication with $A$ is $RG$-isomorphic to $P$. Define the
\emph{Hattori-Stallings rank} of $P$ as
\begin{eqnarray} & \HS_{RG}(P)  =  \tr_{RG}^u(A) &   \in \class_R(G).
\label{Hattori-Stallings_rank}
\end{eqnarray}
The Hattori-Stallings rank depends only on the isomorphism class of
the $RG$-module $P$. It induces an $R$-homomorphism, the
\emph{Hattori-Stallings homomorphism},
\begin{eqnarray}
\HS_{RG} \colon K_0(RG) \otimes_{\IZ} R & \to & \class_R(G), \quad [P] \otimes r~\mapsto~r \cdot \HS_{RG}(P).
\label{HS_K_0(RG)_to_class_0}
\end{eqnarray}

Let $F$ be a field of characteristic zero. Fix an integer $m \ge 1$.
Let $F(\zeta_m) \supset F$ be the Galois extension given by adjoining
the primitive $m$-th root of unity $\zeta_m$ to $F$. Denote by $\Gamma
(m,F)$ the Galois group of this extension of fields, i.e., the group
of automorphisms $\sigma\colon F(\zeta_m) \to F(\zeta_m)$ which induce
the identity on $F$. It can be identified with a subgroup of $\IZ/m^*$
by sending $\sigma$ to the unique element $u({\sigma}) \in \IZ/m^*$
for which $\sigma(\zeta_m) = \zeta_m^{u(\sigma)}$ holds.  Let $g_1$
and $g_2$ be two elements of $G$ of finite order.  We call them
\emph{$F$-conjugate} if for some (and hence all) positive integers $m$
with $g_1^m = g_2^m = 1$ there exists an element $\sigma$ in the
Galois group $\Gamma (m,F)$ with the property that $g_1^{u(\sigma)}$
and $g_2$ are conjugate. Two elements $g_1$ and $g_2$ are
$F$-conjugate for $F = \IQ$, $\IR$ or $\IC$ respectively if the cyclic
subgroups $\langle g_1 \rangle$ and $\langle g_2 \rangle$ are
conjugate, if $g_1$ and $g_2$ or $g_1$ and $g_2^{-1}$ are conjugate,
or if $g_1$ and $g_2$ are conjugate respectively.

Denote by $\con_F(G)_f$ the set of $F$-conjugacy classes $(g)_F$ of
elements $g \in G$ of finite order.  Let $\class_F(G)_f$ be the
$F$-vector space with the set $\con_F(G)_f$ as basis, or,
equivalently, the $F$-vector space of functions $\con_F(G)_f \to F$
with finite support. There are obvious inclusions of $F$-modules
$\class_F(G)_f \subseteq \class_F(G)$.

\begin{lemma} \label{lem:HS_iso_for_finite_G}
  Suppose that $F$ is a field of characteristic zero and $H$ is a
  finite group.  Then the Hattori-Stallings homomorphism (see
  \eqref{HS_K_0(RG)_to_class_0}) induces an isomorphism
  $$\HS_{FH}\colon K_0(FH) \otimes_{\IZ} F \xrightarrow{\cong}
  \class_F(H) = \class_F(H)_f.$$
\end{lemma}
\begin{proof}
  Since $H$ is finite, an $FH$-module is a finitely generated
  projective $FH$-module if and only if it is a (finite-dimensional)
  $H$-representation with coefficients in $F$ and $K_0(FH)$ is the
  same as the representation ring $R_F(H)$.  The Hattori-Stallings
  rank $\HS_{FH}(V)$ and the character $\chi_V$ of a
  $G$-representation $V$ with coefficients in $F$ are related by the
  formula
  $$\chi_V(h)~=~|Z_G\langle h \rangle| \cdot \HS_{FH}(V)(h)$$
  for $h
  \in H$, where $Z_G\langle h \rangle$ is the centralizer of $h$ in
  $G$.  Hence Lemma~\ref{lem:HS_iso_for_finite_G} follows from
  representation theory, see for instance \cite[Corollary 1 on page
  96]{Serre(1997)}.
\end{proof}

Notice that the Bass Conjecture for fields of characteristic zero as
coefficients~\ref{con:Bass_Conjecture_for_FG} is the obvious
generalization of Lemma~\ref{lem:HS_iso_for_finite_G} to infinite
groups.

\begin{lemma} \label{lem:HS_circ_colim_is_injective}
  Suppose that $F$ is a field of characteristic zero and $G$ is a
  group.  Then the composite
  $$\colim_{G/H \in \OrGF{G}{\Fin}} K_0(FH) \otimes_{\IZ} F
  \xrightarrow{I(G,F) \otimes_{\IZ} \id_F} K_0(FG) \otimes_{\IZ} F
  \xrightarrow{\HS_{FG}} \class_F(G)$$
  is injective and has as image
  $\class_F(G)_f$.
\end{lemma}
\begin{proof} This follows from the following commutative diagram, compare~\cite[Lemma~2.15 on page 220]{Lueck(1998b)}.
  $$
  \xymatrix{ \colim_{H \in \OrGF{G}{\Fin}} K_0 (FH ) \otimes_{\IZ}
    F \ar[d]_{\colim_{H \in \OrGF{G}{\Fin}} \HS_{FH}}^{\cong}
    \ar[rr]^{I(G,F) \otimes_{\IZ} \id_F}
    & & K_0 ( F G ) \otimes_{\IZ} F \ar[d]^{\HS_{FG}} \\
    \colim_{H \in \OrGF{G}{\Fin}} \class_F(H)_f \ar[r]^-{j}_-{\cong} &
    \class_F ( G )_f \ar[r]^-{i} & \class_F(G) }$$
  Here the
  isomorphism $j$ is the colimit over the obvious maps $\class_F(H)_f
  \to \class_F(G)_f$ given by extending a class function in the
  trivial way and the map $i$ is the natural inclusion and in
  particular injective.
\end{proof}

\begin{lemma} \label{lem:colim_andH_0_F_field}
  Let $F$ be a field (of arbitrary characteristic).  Then there is a
  natural commutative diagram
  $$\begin{CD} \colim_{\OrGF{G}{\Fin}} K_0(FH) @> \edge > \cong>
    H_0(\EGF{G}{\Fin};K_F)
    \\
    @V I(G,F) VV @VVV
    \\
    K_0(FG) @ > \cong >> H_0(\pt;K_F)
\end{CD}
$$
whose horizontal maps are bijective and whose right vertical arrow
is the assembly map, i.e., the map induced by the projection
$\EGF{G}{\Fin} \to \pt$.
\end{lemma}
\begin{proof}
  This follows from Lemma~\ref{lem:Negative_K-groups_of_FH-Vanish},
  the equivariant Atiyah-Hirzebruch spectral sequence, which converges
  to $H_{p+q}(\EGF{G}{\Fin};K_F)$ in the strong sense and whose
  $E_2$-term $E_{p,q}^2$ is given by the Bredon homology
  $H^{\OrGF{G}{\Fin}}_p(\EGF{G}{\Fin},K_q(F?))$, and the
  natural identification
  $$\colim_{\OrGF{G}{\Fin}}
  K_0(FH)~\cong~H^{\OrGF{G}{\Fin}}_0(\EGF{G}{\Fin},K_0(F?)).$$
\end{proof}

Theorem~\ref{the:FJC_implies_BC_for_FG} follows now from
Lemma~\ref{lem:relative_ass_Fin_to_Vcycv_iso_char_zero},
Lemma~\ref{lem:HS_circ_colim_is_injective} and
Lemma~\ref{lem:colim_andH_0_F_field}.

\begin{lemma} \label{lem:image_of_K_0(RG)_to_K_0(FG)}
  Let $R$ be a commutative integral domain. Let $F$ be its quotient
  field.  Let $G$ be a group.

\begin{enumerate}

\item \label{lem:image_of_K_0(RG)_to_K_0(FG):diagram} Then we obtain a
      commutative diagram
      $$\begin{CD} \colim_{\OrGF{G}{\Fin}} K_0(RH) @>>>
        \colim_{\OrGF{G}{\Fin}} K_0(FH)
        \\
        @V \edge(G,R) VV @V \edge(G,F) VV
        \\
        H_0(\EGF{G}{\Fin};K_R) @>>> H_0(\EGF{G}{\Fin};K_F)
        \\
        @VVV @VVV
        \\
        H_0(\EGF{G}{\VCyc};K_R) @>>> H_0(\EGF{G}{\VCyc};K_F)
        \\
        @VVV @VVV
        \\
        K_0(RG) @>>> K_0(FG)
\end{CD}
$$
where all horizontal maps are change of rings homomorphisms for the
inclusion $R \to F$, the maps $\edge(G,R)$ and $\edge(G,F)$ are edge
homomorphisms appearing in the equivariant Atiyah-Hirzebruch spectral
sequence and the other vertical maps come from the obvious maps
$\EGF{G}{\Fin} \to \EGF{G}{\VCyc}$ and $\EGF{G}{\VCyc} \to \pt$;

\item \label{lem:image_of_K_0(RG)_to_K_0(FG):colim_Fin_R_to_K_0(FG)}
      The image of the composite
      $$\alpha \colon H_0(\EGF{G}{\VCyc};K_R) \otimes_{\IZ} \IQ \to
      H_0(\EGF{G}{\VCyc};K_F) \otimes_{\IZ} \IQ \to K_0(FG)
      \otimes_{\IZ} \IQ$$
      is the same as the image of the composite
      $$\beta \colon \bigoplus_{C \in \FCyc(G)} K_0(RC) \otimes_{\IZ}
      \IQ \to \bigoplus_{C \in \FCyc(G)} K_0(FC) \otimes_{\IZ} \IQ \to
      K_0(FG) \otimes_{\IZ} \IQ$$
      where $\FCyc$ is the family all
      finite cyclic subgroups of $G$.
\end{enumerate}
\end{lemma}

\begin{proof}\ref{lem:image_of_K_0(RG)_to_K_0(FG):diagram} This follows from the
  naturality of the constructions.
  \\[1mm]\ref{lem:image_of_K_0(RG)_to_K_0(FG):colim_Fin_R_to_K_0(FG)}
  For every group $G$,
  every ring $R$, and every $n \in \IZ$ the relative assembly map
\[
H_n^G ( \EGF{G}{\Fin} ; \bfK_R ) \to H_n^G ( \EGF{G}{\VCyc} ; \bfK_R )
\]
is split-injective~\cite{Bartels(2003b)}. This map and the splitting are natural with
respect to change of rings homomorphisms.  Hence
Lemma~\ref{lem:relative_ass_Fin_to_Vcycv_iso_char_zero} implies that the
image of
$$H_n^G ( \EGF{G}{\VCyc} ; \bfK_R ) \otimes_{\IZ} \IQ \to H_n^G (
\EGF{G}{\VCyc} ; \bfK_F ) \otimes_{\IZ} \IQ $$
and the image of the
composite
$$H_n^G ( \EGF{G}{\Fin} ; \bfK_R ) \otimes_{\IZ} \IQ \to H_n^G (
\EGF{G}{\Fin} ; \bfK_F ) \otimes_{\IZ} \IQ \to H_n^G ( \EGF{G}{\VCyc}
; \bfK_F ) \otimes_{\IZ} \IQ$$
agree. By Theorem~1.13 in \cite{Lueck-Reich(2006)}, a variation of the Chern
character in \cite{Lueck(2002b)}, we have 
for every ring $R$ and
group $G$ an isomorphism, natural in $R$,
\begin{multline}
  \chern_n^{G,R} \colon\bigoplus_{p+q = n}~\bigoplus_{\substack{(C)\\C
      \in \FCyc(G)}} H_p(BZ_GC;\IQ) \otimes_{\IQ[W_GC]} \theta_C
  \cdot \left(K_q(RC) \otimes_{\IZ} \IQ \right) \\
  \to~H_n^G(\EGF{G}{\Fin};K_R) \otimes_{\IZ} \IQ,
\label{Chern_character}
\end{multline}
where $W_GC$ is the quotient of the normalizer by the centralizer $Z_GC$ of
the finite cyclic subgroup $C$ in $G$ and $\theta_C$ is an idempotent
in $A(C) \otimes_{\IZ} \IQ$ for $A(C)$ the Burnside ring of $C$ which
singles out a direct summand $\theta_C \cdot
\left(K_q(RC)\otimes_{\IZ} \IQ\right)$ in $K_q(RC) \otimes_{\IZ} \IQ$
by the $A(C) \otimes_{\IZ} \IQ$-module structure on
$K_q(RC)\otimes_{\IZ} \IQ$.  This isomorphism is natural under change
of rings homomorphisms. Recall from
Lemma~\ref{lem:Negative_K-groups_of_FH-Vanish} that $K_q(FC)$ is
trivial for $q \le -1$.  Hence the image of
$$H_0^G ( \EGF{G}{\Fin} ; \bfK_R ) \otimes_{\IZ} \IQ \to H_0^G (
\EGF{G}{\Fin} ; \bfK_F ) \otimes_{\IZ} \IQ$$
agrees with the image of
the composition
\begin{multline*}
  \bigoplus_{(C) \in \FCyc(G)}
  H_0(BZ_GC;\IQ) \otimes_{\IQ[W_GC]} \theta_C \cdot \left(K_0(RC)
    \otimes_{\IZ} \IQ \right)
  \to
  \\
  \bigoplus_{(C) \in \FCyc(G)}
  H_0(BZ_GC;\IQ) \otimes_{\IQ[W_GC]}
  \theta_C \cdot \left(K_0(FC) \otimes_{\IZ} \IQ \right)
  \\
  = \bigoplus_{p+q = 0} \bigoplus_{(C) \in \FCyc(G)} H_p(BZ_GC;\IQ)
  \otimes_{\IQ[W_GC]} \theta_C \cdot \left(K_q(FC) \otimes_{\IZ} \IQ
  \right)
  \\
  \xrightarrow{\chern^{G,F}_0} H_0^G(\EGF{G}{\Fin};K_F) \otimes_{\IZ}
  \IQ.
\end{multline*}

Since $H_0(BZ_GC;\IQ)$ is $\IQ$ with the trivial $W_GC$-action,
$H_0(BZ_GC;\IQ) \otimes_{\IQ[W_GC]} \theta_C \cdot
\left(K_0(RC)\otimes_{\IZ} \IQ \right)$ is a quotient of $\theta_C
\cdot \left(K_0(RC)\otimes_{\IZ} \IQ \right)$ which is a direct
summand in $K_0(RC)\otimes_{\IZ} \IQ$.  Hence $\im(\alpha) \subseteq
\im(\beta)$. We get $\im(\beta) \subseteq \im(\alpha)$ from the
commutative diagram appearing in
assertion~\ref{lem:image_of_K_0(RG)_to_K_0(FG):diagram}. This finishes
the proof of Lemma~\ref{lem:image_of_K_0(RG)_to_K_0(FG)}
\end{proof}

Finally we can give the proof of Theorem~\ref{the:Bass_for_R_to_F}.
\begin{proof}
  For any finite group of $H \subseteq G$ such that its order is not invertible in $R$ the
  map $\widetilde{K}_0(RH) \to \widetilde{K}_0(FH)$ is trivial by a result of
  Swan~\cite[Theorem~8.1]{Swan(1960a)} (see also~\cite[Corollary~4.2]{Bass(1979)}).
  Lemma~\ref{lem:image_of_K_0(RG)_to_K_0(FG)}~\ref{lem:image_of_K_0(RG)_to_K_0(FG):colim_Fin_R_to_K_0(FG)}
  implies that the composite
  $$H_0(\EGF{G}{\VCyc};K_R) \otimes_{\IZ} \IQ \to H_0(\EGF{G}{\VCyc};K_F) \otimes_{\IZ}
  \IQ \to K_0(FG) \otimes_{\IZ} \IQ \to \widetilde{K}_0(FG) \otimes_{\IZ} \IQ$$
  is trivial.
  Since $G$ satisfies the Farrell-Jones Conjecture for algebraic $K$-theory with
  coefficients in $R$ by assumption, the map
  $H_0(\EGF{G}{\VCyc};K_R) \otimes_{\IZ} \IQ~\to~K_0(RG) \otimes_{\IZ} \IQ$
  is   surjective. Because of the commutative diagram appearing in
  Lemma~\ref{lem:image_of_K_0(RG)_to_K_0(FG)}~\ref{lem:image_of_K_0(RG)_to_K_0(FG):diagram}
  the map
  $$K_0(RG)\otimes_{\IZ} \IQ \to K_0(FG)\otimes_{\IZ} \IQ \to
  \widetilde{K}_0(FG) \otimes_{\IZ} \IQ$$
  is trivial. Now Theorem~\ref{the:Bass_for_R_to_F}  follows.
\end{proof}


\section{The Kaplansky Conjecture}
\label{sec:The_Kaplansky_Conjecture}

We need the following definition.

\begin{definition}\label{def:directly_finite}
  An $R$-module $M$ is called \emph{directly finite}  if every
  $R$-module $N$ satisfying $M \cong_R M \oplus N$ is trivial. A ring
  $R$ is called \emph{directly finite} (or \emph{von Neumann finite})
  if it is directly finite as a
  module over itself, i.e., if $r,s \in R$ satisfy $rs= 1$, then $sr =
  1$. A ring is called \emph{stably finite} if the matrix algebra $M(n,n,R)$
  is directly finite for all $n \ge 1$.
\end{definition}

\begin{remark} \label{rem:stable-finite-detects-zero}
Stable finiteness for a ring  $R$ is equivalent to the following statement.
Every finitely generated projective $R$-module $P$ 
whose class in $K_0(R)$ is zero is already the trivial module, i.e., $0 = [P] \in K_0 ( R)$ implies $P \cong 0$.
\end{remark}

\begin{theorem}\label{the:examples_for_stably_finite_group_rings}

\begin{enumerate}

\item \label{the:examples_for_stably_finite_group_rings:fields_of_characteristic_zero}
If $F$ is a field of characteristic zero, then $FG$ is  stably finite for every group $G$;

\item \label{the:examples_for_stably_finite_group_rings:sofic_groups}
If $R$ is a skew-field and $G$ is a sofic group, then $RG$ is stably finite.

\end{enumerate}
\end{theorem}
\begin{proof}\ref{the:examples_for_stably_finite_group_rings:fields_of_characteristic_zero}
This is proved by Kaplansky~\cite{Kaplanski(1969)} (see also
Passman~\cite[Corollary~1.9 on page 38]{Passman(1977)}).
\\[1mm]\ref{the:examples_for_stably_finite_group_rings:sofic_groups}
This is proved for free-by-amenable groups by Ara-Meara-Perera~\cite{Ara-Meara-Perera(2002)} and extended
to sofic groups by Elek-Szabo~\cite[Corollary~4.7]{Elek-Szabo(2004)}.
\end{proof}

\begin{lemma} \label{lem:widetildeK_0(RG)_is_0_and_Kaplansky}
  Let $R$ be a ring whose idempotents are all trivial.  Let $G$ be a
  group such that the map induced by the inclusion $i \colon R \to RG$
  $$i_* \colon K_0(R) \otimes_{\IZ} \IQ \xrightarrow{} K_0(RG) \otimes_{\IZ} \IQ$$
  is bijective and $RG$ is stably finite.

  Then all idempotents in $RG$ are trivial.
\end{lemma}

\begin{proof}
  Let $p$ be an idempotent in $RG$. We want to show $p \in \{0,1\}$.
  Denote by $\epsilon \colon RG \to R$ the augmentation homomorphism
  sending $\sum_{g \in G} r_g \cdot g$ to $\sum_{g \in G} r_g$.  Since
  $\epsilon(p) \in R$ is $0$ or $1$ by assumption, it suffices to show
  $p = 0$ under the assumption that $\epsilon(p) = 0$.  Let $(p)
  \subseteq RG$ be the ideal generated by $p$ which is a finitely
  generated projective $RG$-module. Since $i_* \colon K_0(R)
  \otimes_{\IZ} \IQ \to K_0(RG) \otimes_{\IZ} \IQ$ is surjective
  by assumption, we can find a finitely generated projective
  $R$-module $P$ and integers $k,m,n \ge 0$ satisfying
  $$(p)^k \oplus RG^m \cong_{RG} i_*(P) \oplus RG^n.$$
  If we now apply $i_{\ast} \circ \epsilon_{\ast}$ and use 
  $\epsilon \circ i = \id$, $i_{\ast} \circ \epsilon_{\ast} ( RG^l ) \cong RG^l$ and $\epsilon(p) =0$
  we obtain
  $$RG^m \cong i_{\ast} (P) \oplus RG^n.$$
  Inserting this in the first equation yields
  $$(p)^k \oplus i_{\ast}( P ) \oplus RG^n \cong i_{\ast} ( P ) \oplus RG^n$$ and therefore
  $0 = [ (p)^k ] \in K_0(RG)$. Using Remark~\ref{rem:stable-finite-detects-zero} we conclude
  $(p)^k = 0$ and hence $p = 0$.
\end{proof}

%

Now Theorem~\ref{the:FJC_and_Kaplansky} follows from
Theorem~\ref{the:examples_for_stably_finite_group_rings} and
Lemma~\ref{lem:widetildeK_0(RG)_is_0_and_Kaplansky}.


\section{Nil-groups}
\label{sec:Nil-groups}

In this section we give a brief discussion about Waldhausen's Nil-groups.

\subsection{Applications of homotopy $K$-theory to Waldhausen's Nil-groups}
\label{subsec:applications-KH-to-Waldhausen-Nil}
The Nil-groups due to Bass $N\!K_n(R)$ have been generalized to
\emph{Waldhausen's Nil groups} as follows.  These groups were defined by
Waldhausen (see~\cite{Waldhausen(1978a)}, \cite{Waldhausen(1978b)})
for $n \ge 0$ and their extensions to $n \in \IZ$ is explained
in~\cite[Section~9]{Bartels-Lueck(2004ind)}. For more details we
refer to these papers. There are two kinds of Nil-groups,
the first one taking care of amalgamated products and the second of generalized
Laurent extensions.

We begin with a discussion of the one for amalgamated products. Let
$\alpha \colon C \to A$ and $\beta\colon C \to B$ be pure and free and
write $A = \alpha(C) \oplus A'$ and $B = \beta(C) \oplus B'$ and put
$R = A \ast_C B$. Then the Nil-groups of the first kind are denoted by
$\Nil_n(C;A',B')$. The group $\Nil_n(C;A',B')$ is a direct summand in
$K_{n+1}(R)$ and there is a long exact sequence
\begin{multline}
  \cdots \to K_{n+1}(RC) \to K_{n+1}(A) \oplus K_{n+1}(B) \to
  K_{n+1}(R)/\Nil_n(C;A',B')
  \\
  \to K_n(C) \to K_{n}(A) \oplus K_{n}(B) \to \cdots.
  \label{sequence_for_K(A_ast_C_B)}
\end{multline}
There is an analogous long exact sequence for homotopy $K$-theory which does not
involve Nil-terms~(see \cite[Section~9]{Bartels-Lueck(2004ind)})
\begin{multline}
  \cdots \to \KH_{n+1}(RC) \to \KH_{n+1}(A) \oplus \KH_{n+1}(B) \to
  \KH_{n+1}(R)
  \\
  \to \KH_n(C) \to \KH_{n}(A) \oplus \KH_{n}(B) \to \cdots.
   \label{sequence_for_KH(A_ast_C_B)}
\end{multline}
The natural transformation from algebraic $K$-theory to homotopy
$K$-theory induces a map between the long exact
sequences~\eqref{sequence_for_K(A_ast_C_B)}
and~\eqref{sequence_for_KH(A_ast_C_B)}. Now suppose that for $\IZ
\subseteq \Lambda \subseteq \IQ$ the maps
\begin{eqnarray*}
K_n(A) \otimes_{\IZ} \Lambda & \to & \KH_n(A) \otimes_{\IZ} \Lambda;\\
K_n(B) \otimes_{\IZ} \Lambda & \to & \KH_n(B) \otimes_{\IZ} \Lambda;\\
K_n(C) \otimes_{\IZ} \Lambda & \to & \KH_n(C) \otimes_{\IZ} \Lambda;\\
K_n(R) \otimes_{\IZ} \Lambda & \to & \KH_n(R) \otimes_{\IZ} \Lambda,
\end{eqnarray*}
are bijective for all $n \in \IZ$. Then a Five-Lemma argument implies
$$\Nil(A;B',C') \otimes_{\IZ} \Lambda = 0 \quad \text{ for all } n \in \IZ.$$
Recall from the spectral sequence of~\cite[Theorem~1.3]{Weibel(1989)}
that $K_n(A) \otimes_{\IZ} \Lambda \to \KH_n(A) \otimes_{\IZ} \Lambda$
is bijective if we have $N\!K_n(A[x_1,x_2, \ldots ,x_k])
\otimes_{\IZ} \Lambda = 0$ for all $n \in \IZ$ and $k \in \IZ, k \ge 0$ and analogous for $B$,
$C$ and $R$ instead of $A$.

Next we deal with generalized Laurent extensions. Let $R$ be the
generalized Laurent extension of pure and free ring maps
$\alpha,\beta \colon C \to A$.  Associated to it is the  Nil-term
$\Nil_n(C;{_{\alpha}A'_{\alpha}},{_{\beta}A''_{\beta}},{_{\beta}A_{\alpha}},{_{\alpha}A_{\beta}})$
which is a  direct summand in $K_{n+1}(R)$. We obtain long exact sequences
\begin{multline}
  \cdots \to K_{n+1}(C) \xrightarrow{\alpha_* - \beta_*} K_{n+1}(A)
  \to
  K_{n+1}(R)/\Nil_n(C;{_{\alpha}A'_{\alpha}},{_{\beta}A''_{\beta}},{_{\beta}A_{\alpha}},{_{\alpha}A_{\beta}})
  \\
  \to K_{n}(C) \xrightarrow{\alpha_* - \beta_*} K_{n}(A) \to
  \cdots \label{sequence_for_K(Laurent)}
\end{multline}
and
\begin{multline}
  \cdots \to \KH_{n+1}(C) \xrightarrow{\alpha_* - \beta_*} \KH_{n+1}(A)
  \to \KH_{n+1}(R)
  \\
  \to \KH_{n}(C) \xrightarrow{\alpha_* - \beta_*} \KH_{n}(A) \to
  \cdots.  \label{sequence_for_KH(Laurent)}
\end{multline}
If for $\IZ \subseteq \Lambda \subseteq \IQ$ the maps
\begin{eqnarray*}
K_n(A) \otimes_{\IZ} \Lambda & \to & \KH_n(A) \otimes_{\IZ} \Lambda;\\
K_n(C) \otimes_{\IZ} \Lambda & \to & \KH_n(C) \otimes_{\IZ} \Lambda;\\
K_n(R) \otimes_{\IZ} \Lambda & \to & \KH_n(R) \otimes_{\IZ} \Lambda,
\end{eqnarray*}
are bijective for all $n \in \IZ$,
then we conclude as above that
$$\Nil_n(C;{_{\alpha}A'_{\alpha}},{_{\beta}A''_{\beta}},{_{\beta}A_{\alpha}},{_{\alpha}A_{\beta}})
\otimes_{\IZ} \Lambda = 0 \quad \text{ for all } n \in \IZ$$
holds.

We can now prove Theorem~\ref{thm:waldhausen-nil-is-rationally-trivial}.

\begin{proof}\ref{thm:waldhausen-nil:amalgamated}.
  Because $G$, $H$ and $K$ are finite we conclude from Lemma~\ref{lem:K-to-KH}~\ref{lem:K-to-KH:IZ}
  that the maps
  \begin{eqnarray*}
  K_n(\IZ G) \otimes_\IZ \IQ & \to & \KH_n(\IZ G) \otimes_\IZ \IQ; \\
  K_n(\IZ H) \otimes_\IZ \IQ & \to & \KH_n(\IZ H) \otimes_\IZ \IQ; \\
  K_n(\IZ K) \otimes_\IZ \IQ & \to & \KH_n(\IZ K) \otimes_\IZ \IQ,
  \end{eqnarray*}
  are bijective.
  Let $\Gamma := G *_K H$ be the amalgamated product.
  Then $\Gamma$ acts properly and cocompactly on a tree and is therefore word-hyperbolic.
  Theorems~\ref{the:class_FJC}~\ref{the:class_FJC:word_hyper_virt_nilpotent}
  and~\ref{the:class_FJ_KH_andFJ_N}~\ref{the:class_FJ_KH_andFJ_N:elem_amen_word_hyper}
  imply $\Gamma \in \calfj(\IZ)$ and $\Gamma \in \calfj_{\KH}(\IZ)$.
  Consider the following  commutative diagram
  $$
  \xycomsquare{H_n^\Gamma(\EGF{\Gamma}{\Fin};\bfK_{\IZ})\otimes_\IZ \IQ}{}{H_n(\pt;\bfK_{\IZ})\otimes_\IZ \IQ
    = K_n(R\Gamma) \otimes_\IZ \IQ} {}{}
  {H_n^\Gamma(\EGF{\Gamma}{\Fin};\bfKH_{\IZ})\otimes_\IZ \IQ}{}{H_n(\pt;\bfKH_{\IZ})\otimes_\IZ \IQ =
    \KH_n(R\Gamma) \otimes_\IZ \IQ}
  $$
  where the horizontal maps are the assembly maps induced by the
  projection $\EGF{\Gamma}{\Fin} \to \pt$ and the vertical maps are
  induced by the passage from algebraic $K$-theory to homotopy $K$-theory.
  The lower horizontal map is an isomorphism since $\Gamma \in \calfj_{\KH}(\IZ)$,
  the upper horizontal map is an isomorphism by Lemma~\ref{lem:relative_ass_Fin_to_Vcycv_R_is_Z}
  since $\Gamma \in \calfj(\IZ)$.
  The left vertical map is an isomorphism by
  Lemma~\ref{lem:assembly-for-K-versus-KH}~\ref{lem:assembly-for-K-versus-KH:IZ}.
  Therefore $K_n(\IZ \Gamma) \otimes_{\IZ} \IQ \cong \KH_n(\IZ \Gamma) \otimes_{\IZ} \IQ$.
  The Five-Lemma argument from above implies now the rational vanishing of the Nil-groups as claimed.
\\[1mm]\ref{thm:waldhausen-nil:laurent}
  Here we can argue as in~\ref{thm:waldhausen-nil:amalgamated} but now using
  the $HNN$-extension associated to the inclusions $\alpha$ and $\beta$
  instead of the amalgamated product.
\end{proof}


\subsection{Nil-groups and rings with finite characteristic}
\label{subsec:Nil-groups_and_rings_with_finite_characteristic}

For this subsection fix an integer $N \ge 2$  and a ring $A$
of characteristic $N$.  In
Lemma~\ref{lem:Vanishing_of_NK_p(R)[1/N]} we have shown that we get
$N\!K_n(A)[1/N]=0$ for every $n \in \IZ$.
This implies by the above Five-Lemma argument for all $n \in \IZ$ and all pure and free maps $\alpha
\colon C \to A$ and $\beta\colon C \to B$ and all pure and free maps
$\alpha,\beta \colon C \to A$ respectively
\begin{eqnarray*}
\Nil_n(C;A',B') \otimes_{\IZ}[1/N] & = & 0;
\\
\Nil_n(C;{_{\alpha}A'_{\alpha}},{_{\beta}A''_{\beta}},{_{\beta}A_{\alpha}},{_{\alpha}A_{\beta}})
\otimes_{\IZ}[1/N] & = & 0.
\end{eqnarray*}


\subsection{Nil-groups and regular rings $R$ with $\IQ \subset R$}
\label{subsec:Nil-groups_and_regular_rings_R_with_Z_subseteq_R}

\begin{conjecture} \label{con:NK_n(RG)_for_regular_R}
Let $R$ be a regular ring with $\IQ \subseteq R$. Then we get
for all groups $G$ and all $n \in \IZ$ 
$$N\!K_n(RG) = 0$$
and  the canonical map
$$K_n(RG) \to \KH_n(RG)$$
is bijective.
\end{conjecture}

By the following discussion
Conjecture~\ref{con:NK_n(RG)_for_regular_R}  is  true if $G$ belongs to the class
$\calfj(R)$ appearing in Theorem~\ref{the:class_FJC} but we do not
know it for all groups $G$.

If $R$ is regular then the polynomial ring $R[t]$ is regular as well.
Thus if $N\!K_n(RG) = 0$ holds for a fixed group $G$ and all regular rings $R$ with $\IQ \subseteq R$,
then $K_n(RG) \to \KH_n(RG)$ is bijective for this fixed group $G$ and all regular rings $R$ with $\IQ \subseteq R$.
This follows from the spectral sequence appearing in~\cite[Theorem~1.3]{Weibel(1989)}.
Using the Bass-Heller-Swan decomposition it is not hard to see that similarly,
if $N\!K_n(RG) = 0$ holds for all groups $G$ and a fixed regular ring $R$ with $\IQ \subseteq R$,
then $K_n(RG) \to \KH_n(RG)$ is bijective for all groups $G$ and this fixed regular ring $R$ with $\IQ \subseteq R$.

If Conjecture~\ref{con:NK_n(RG)_for_regular_R} is true, then also Waldhausen's Nil-groups
associated to a free amalgamated product of groups $G = G_1 \ast_{G_0} G_1$ or to an
HNN-extension with coefficients in a regular ring $R$ with $\IQ \subseteq R$ vanish by
the general argument from Section~\ref{subsec:applications-KH-to-Waldhausen-Nil}.

Let $\calfj_{Fin}(R)$ be the class of those
groups for which the Fibered Farrell-Jones Conjecture for algebraic $K$-theory with
coefficients in $R$ is true for the family $\Fin$. We get $\calfj_{Fin}(R) \subseteq
\calfj(R)$ from Lemma~\ref{lem:enlarging_the_family_for_the_Fibered_Isomorphism_Conjecture}
for all rings $R$ but $\calfj_{Fin}(\IZ) \not= \calfj(\IZ)$.

If Conjecture~\ref{con:NK_n(RG)_for_regular_R} holds and $R$ is a regular ring with $\IQ \subseteq R$,
then the map $H^G_n( X; \bfK_R) \to H^G_n( X; \bfKH_R)$ induced by the transition
from $K$-theory to homotopy $K$-theory is an isomorphism for
any $G$-$CW$-complex $X$ and the assembly maps for $K$-theory and homotopy $K$-theory can be
identified.
Thus Conjecture~\ref{con:NK_n(RG)_for_regular_R} implies $\calfj_{Fin}(R) = \calfj_{\KH}(R)$.
In particular, if  Conjecture~\ref{con:NK_n(RG)_for_regular_R} holds,
then by Theorem~\ref{the:class_FJ_KH_andFJ_N}~\ref{the:class_FJ_KH_andFJ_N:actions_on_trees}
every virtually cyclic group satisfies the Fibered Farrell-Jones Conjecture for algebraic $K$-theory with
coefficients in $R$ for the family $\Fin$.
We omit the proof (based on the Bass-Heller-Swan formula) that
Conjecture~\ref{con:NK_n(RG)_for_regular_R} is in fact equivalent
to the statement $\IZ \in \calfj_{Fin}(R)$ for any regular ring $R$ with $\IQ \subseteq R$.

If Conjecture~\ref{con:NK_n(RG)_for_regular_R} is true and $R$ is a regular ring with $\IQ \subseteq R$,
then the Transitivity Principle (see Theorem~\ref{the:transitivity})
implies that
$\calfj(R) = \calfj_{Fin}(R)$, the conclusions appearing in Theorem~\ref{the:class_FJ_KH_andFJ_N} are
true for each group in $\calfj(R) = \calfj_{Fin}(R)$, and
assertion~\ref{the:FJC_implies_Moodys_induction:reg_ring_Q_subset_R} of
Theorem~\ref{the:FJC_implies_Moodys_induction} is true for all groups in $\calfj_{Fin}(R)$
and in particular for all elementary amenable groups.

Let $R$ be any ring. Recall that  $\calfj(R)$ contains for
instance word-hyperbolic groups and virtually nilpotent groups by
Theorem~\ref{the:class_FJC} but we do not know whether elementary amenable groups belong
to $\calfj(R)$ and whether they satisfy
assertion~\ref{the:FJC_implies_Moodys_induction:reg_ring_Q_subset_R} of
Theorem~\ref{the:FJC_implies_Moodys_induction}. Our methods do not give this conclusion
since we know only that $G \in \calfj(R)$ implies $G \times \IZ \in \calfj(R)$ but we
would at least need $G \rtimes_{\alpha} \IZ \in \calfj(R)$ for any automorphism $\alpha
\colon G \to G$. More generally, it would be interesting to solve the question whether for
any extensions $ 1 \to G \to \Gamma \to V \to 1$ for virtually cyclic $V$ and $G \in
\calfj(R)$ we have $\Gamma \in \calfj(R)$ because a positive answer would imply
that $\calfj(R)$ is closed under extensions.


\section{The Farrell-Jones Conjecture for $L$-theory}
\label{sec:The_Farrell-Jones_Conjecture_for_L-theory}

In this section ring $R$  will always mean associative ring with unit and involution.

Recall that the \emph{(Fibered) Farrell-Jones Conjecture for algebraic
  $L$-theory with coefficients in $R$ for the group $G$} is the
(Fibered) Isomorphism
Conjecture~\ref{def:(Fibered)_Isomorphism_Conjectures_for_calh?_ast}
in the special case, where the family $\calf$ consists of all
virtually cyclic subgroups of $G$ and $\calh^?_*$ is the equivariant
homology theory $\calh^?_*(-;\bfL_R^{\langle - \infty \rangle})$
associated to the $\Groupoids$-spectrum given by algebraic $L$-theory
and $R$ as coefficient ring (see \cite[Section~6]{Lueck-Reich(2005)}).
So the Farrell-Jones Conjecture for algebraic $L$-theory with
coefficients in $R$ for the group $G$ predicts that the map
$$H_n^G(\EGF{G}{\VCyc},\bfL_R^{\langle - \infty \rangle}) \to
L_n^{\langle -\infty \rangle}(RG)$$
is bijective for all $n \in \IZ$.
The original source for (Fibered) Farrell-Jones Conjecture is
\cite{Farrell-Jones(1993a)}. The corresponding conjecture is false if
one replaces the decoration $\langle - \infty \rangle$ with the
decoration $p$, $h$ or $s$ (see~\cite{Farrell-Jones-Lueck(2002)}).
For the status of the Farrell-Jones Conjecture with coefficients in
$\IZ$  we refer for instance to \cite[Sections~5.2 and~5.3]{Lueck-Reich(2005)}.

The next result is proved in~\cite[Lemma~5.2]{Bartels-Echterhoff-Lueck(2007colim)}.
\begin{lemma} \label{lem:H(-,L_R-infty)_strongly_continuous)} Let $R$ be a ring.
  Then the equivariant homology theory $\calh^?_*(-;\bfL_R^{\langle -\infty\rangle})$ is strongly
  continuous.
\end{lemma}

Theorem~\ref{the:isomorphism_conjecture_is_stable_under_colim}~\ref{the:isomorphism_conjecture_is_stable_under_colim:general}
and
Lemma~\ref{lem:FJC_for_L[1/2]_Fin_to_Vcyc} imply
that for any ring $R$ and any direct system of groups $\{G_i
\mid i \in I\}$ (with not necessarily injective structure maps)
$G = \colim_{i \in I} G_i$ satisfies the Fibered Farrell Jones Conjecture for
algebraic $L$-theory with coefficients in $R$ if each group $G_i$
does.
 
In algebraic $L$-theory the situation simplifies if one inverts $2$.

\begin{lemma} \label{lem:FJC_for_L[1/2]_Fin_to_Vcyc}
Let $R$ be a ring and let $G$ be a group. Then
the Fibered Farrell Jones Conjecture with coefficients in $R$
after applying $- \otimes_\IZ \IZ[1/2]$ to the assembly map
holds for $(G,\VCyc)$
if and only the same holds for
$(G,\Fin)$.
\end{lemma}

\begin{proof}
Because of  Lemma~\ref{lem:extensions_with_finite_kernel}
and the fact that every virtually cyclic group maps surjectively onto $\IZ$ or $\IZ/2 \ast \IZ/2$
it suffices to show the claim for $G = \IZ$ and $G = \IZ/2 \ast \IZ/2$.
These cases follow from the exact sequence involving UNil-terms and the proof that UNil-groups
are $2$-torsion due to Cappell~\cite{Cappell(1974c)} in the case
$\IZ/2 \ast \IZ/2$. In the case $\IZ$ there are no  UNil-terms for infinite virtually cyclic
groups of the first kind.  This follows essentially
from~\cite{Ranicki(1973b)} and~\cite{Ranicki(1973c)}
as carried out in~\cite[Lemma~4.2]{Lueck(2005heis)}.
\end{proof}

The Transitivity Principle (see Theorem~\ref{the:transitivity}) implies
that for a ring $R$ and a group $G$ the Farrell-Jones Conjecture for algebraic
$L$-theory with coefficients in $R$ is true after applying $-
\otimes_{\IZ} \IZ[1/2]$ to the assembly map for the family $\VCyc$ if
and only if this is true for the family $\Fin$. This implies that
assertions~\ref{the:class_FJ_KH_andFJ_N:products},\ref{the:class_FJ_KH_andFJ_N:colim},%
\ref{the:class_FJ_KH_andFJ_N:subgroups} and~\ref{the:class_FJ_KH_andFJ_N:actions_on_trees} appearing in
Theorem~\ref{the:class_FJ_KH_andFJ_N} are true for the family
$\calfj_{L[1/2]}(R)$ of those groups for which the Farrell-Jones
Conjecture for algebraic $L$-theory with coefficients in $R$ is true
after applying $- \otimes_{\IZ} \IZ[1/2]$. We do not know whether all
word hyperbolic groups belong to $\calfj_{L[1/2]}(R)$. Farrell and
Jones claim in \cite[Remark~2.1.3]{Farrell-Jones(1993a)} without gi©ving 
details that the $L$-theory version of their
celebrated~\cite[Theorem~1.2]{Farrell-Jones(1993a)} is true.  This
together with \cite[3.6.4]{Ranicki(1981)} implies that
all virtually finitely generated abelian groups belong to
$\calfj_{L[1/2]}(R)$ for $R = \IZ, \IQ$.  Hence all elementary amenable groups belong to
$\calfj_{L[1/2]}(R)$ for $R = \IZ,\IQ$ by
Lemma~\ref{lem:elementary_amenable}.  This has already been proved
in~\cite[Section~5]{Farrell-Linnell(2003a)}.

\def\cprime{$'$} \def\polhk#1{\setbox0=\hbox{#1}{\ooalign{\hidewidth
  \lower1.5ex\hbox{`}\hidewidth\crcr\unhbox0}}}

 \end{document}